\documentclass[11pt]{amsart}

\usepackage{amstext,amssymb,amsmath,amsbsy}

\usepackage{amsmath,amssymb,latexsym,dsfont}
\usepackage[left=2.5cm,right=2.5cm,top=2.2cm,bottom=2.2cm]{geometry}

\usepackage{amsmath,amsfonts,amssymb,epsfig, textcomp}
\usepackage{graphicx,tikz}
\usepackage{hyperref, cleveref} 
\usepackage{cases,listings}
\usepackage{color}
\usepackage{mathabx}
\usepackage{esint}

\newtheorem{theorem}{Theorem}
\newtheorem{lemma}{Lemma}
\newtheorem{conjecture}{Conjecture}

\newtheorem{proposition}{Proposition}

\newtheorem{remark}{Remark}

\newtheorem{question}{Open question}

\theoremstyle{definition}
\numberwithin{equation}{section}
\numberwithin{theorem}{section}
\numberwithin{lemma}{section}
\numberwithin{proposition}{section}
\numberwithin{remark}{section}
\numberwithin{example}{section}
\numberwithin{corollary}{section}
\numberwithin{exercise}{section}
\numberwithin{definition}{section}

\newcommand{\card}[1]{\mathrm{card}(#1)}
\newcommand{\lr}{\rightarrow}
\newcommand{\dsp}{\displaystyle}
\newcommand{\esssup}{\dsp \mathop{\mathrm{ess} \sup}_{ I}}
\newcommand{\essinf}{\dsp \mathop{\mathrm{ess} \inf}_{ I}}

\newcommand{\sge}{\gtrsim}

\newcommand{\supp}{\operatorname{supp}}

\newcommand{\dive}{\operatorname{div}}

\newcommand{\eps}{\varepsilon}

\newcommand{\loc}{_{loc}}

\newcommand{\mN}{\mathbb{N}}

\newcommand{\mS}{\mathbb{S}}
\newcommand{\mR}{\mathbb{R}}

\newcommand{\ZZ}{\mathbb{Z}}
\newcommand{\mZ}{\mathbb{Z}}

\newcommand{\cM}{{\mathcal M}}

\newcommand{\hmr}[1]{\textcolor{black}{#1}}

\newcommand{\be}{\begin{equation}}
\newcommand{\ee}{\end{equation}}

\newcommand{\cQ}{{\mathcal Q}}

\newcommand{\ka}{\kappa}

\title[Characterizations of the Sobolev norms and the total variation]{Characterizations of the Sobolev norms and the total variation via nonlocal functionals, and related problems}

\author[H.-M. Nguyen]{Hoai-Minh Nguyen}
\address[H.-M. Nguyen]{Sorbonne Universit\'e, Universits\'e Paris Cit\'e, CNRS, INRIA, \newline
\indent Laboratoire Jacques-Louis Lions, LJLL, F-75005 Paris, France
}
\email{hoai-minh.nguyen@sorbonne-universite.fr}

\begin{document}

\maketitle

\begin{center}
    \textit{Dedicated to Ha\"{\i}m Brezis with admiration, gratitude, and memories}
\end{center}

\begin{abstract}
We  briefly discuss the contribution of Ha\"im Brezis and his co-authors on the characterizations of the Sobolev norms and the total variation using non-local functionals. Some ideas of the analysis are  given and new results are presented. 
\end{abstract}



\tableofcontents

\section{Introduction} \label{sect-introduction}

In this paper, we first briefly discuss the contribution of Ha\"im Brezis and his co-authors on the characterizations of the Sobolev norms and the total variation using non-local functionals, and related problems. Some ideas of the analysis are also given. We then present several new results by developing these ideas. 

\medskip 
We begin with the BBM formula due to Bourgain, Brezis, and Mironescu \cite{BBM-01} (see also \cite{Brezis-02, Davila-02}). To this end, for $N \ge 1$, $p \ge 1$,  and $u \in L^p(\mR^N)$, it is convenient to denote
\be
\Phi (u) =  \left\{  \begin{array}{c} \|\nabla u \|_{L^p(\mR^N)}^p \, dx \mbox{ if } p > 1,\\[6pt]
 \| \nabla u \|_{\cM(\mR)} \mbox{ if } p =1.   
\end{array}\right. 
\ee
Recall that, for $f \in L^1(\mR^N)$, 
$$
\|\nabla f \|_{\cM(\mR^N)} : = \sup \left\{ \left| \int_{\mR^N} f \dive \varphi \right| ; \varphi \in C^\infty_c(\mR^N) \mbox{ with } \| \varphi \|_{L^\infty(\mR^N)}  \le 1\right\}. 
$$
In what follows, a sequence of functions  $(\rho_n)_{n
\ge 1} \subset L^1(0, + \infty)$ is called a sequence of non-negative mollifiers if the following properties hold: 
\begin{equation*}
\begin{array}{ccc}
\rho_n \geq 0, 
\end{array}
\end{equation*}
\begin{equation*}
\dsp \lim_{ n \lr \infty } \int_\tau^\infty \rho_n(r) r^{N-1} \,
dr = 0 \quad   \forall \, \tau >0, \quad \mbox{ and } \quad
 \int_{0}^{+ \infty} \rho_n(r) r^{N-1} \, dr =
1.
\end{equation*}

Here is the BBM formula.  

\begin{theorem}[BBM formula, Bourgain \& Brezis \& Mironescu]\label{thm-BBM}
Let $N \ge 1$, $1 \le  p < + \infty$, and let $(\rho_n)_{n \ge 1}$ be a sequence of non-negative  mollifiers.  Then, for $u \in L^p(\mR^N)$,  
\be\label{thm-BBM-cl1}
\mathop{\iint}_{\mR^N \times \mR^N} \frac{|u(x)
- u(y)|^p}{|x-y|^p} \rho_n(|x-y|)\, dx \, dy \le  C_{N, p} \Phi (u), 
\ee
and 
\begin{equation}\label{thm-BBM-cl2}
\lim_{n  \lr  \infty } \mathop{\iint}_{\mR^N \times \mR^N} \frac{|u(x)
-u(y)|^p}{|x-y|^p} \rho_n(|x-y|)\, dx \, dy = K_{N, p} \Phi(u), 
\end{equation}
where $K_{N,p}$ is defined by 
\begin{equation}\label{def-K}
K_{N,p}= \int_{\mS^{N-1}} |e \cdot \sigma|^p \, d \sigma,
\end{equation}
for any $e \in \mS^{N-1}$, the unit sphere of $\mR^N$. 
\end{theorem}

Here and in what follows, $C_{N, p}$ denotes a positive constant depending only on $N$ and $p$, and might change from one place to another. 

\medskip

As a convention, the RHS of \eqref{thm-BBM-cl1} or \eqref{thm-BBM-cl2} is infinite  if $u\not \in W^{1, p}(\mR^N)$ for $p> 1$ and $u \not \in BV(\mR^N)$ for $p=1$. 

\medskip 
\Cref{thm-BBM} was established by Bourgain, Brezis, and Mironesu \cite{BBM-01} (see also \cite{Brezis-02}) in the case $p>1$. In the case $p=1$, they also showed there that the $\liminf$ and the $\limsup$ of the LHS of \eqref{thm-BBM-cl1} as $n \to + \infty$ is of the order of the RHS of \eqref{thm-BBM-cl1} instead of \eqref{thm-BBM-cl2}. The proof of \eqref{thm-BBM-cl2} in the case $p=1$ and $u \in BV(\mR^N)$ is due to Davila \cite{Davila-02}. 

\medskip 
We next briefly discuss some ideas of the proof. We \hmr{first} deal with \eqref{thm-BBM-cl1} under the additional assumption that $\Phi(u) < + \infty$. We first consider the case where $u \in C^\infty_c(\mR^N)$. Using the fact, for $x, y \in \mR^N$,  
$$
u(y) - u(x) = \int_0^1 \nabla u(x + t (y -x)) \cdot (y-x) \, dt, 
$$
and Jensen's inequality, one can prove that 
\begin{multline} \label{thm-BBM-p1}
\mathop{\iint}_{\mR^N \times \mR^N} \frac{|u(x)
- u(y)|^p}{|x-y|^p} \rho_n(|x-y|)\, dx \, dy \\[6pt] 
\le \mathop{\iint}_{\mR^N \times \mR^N}  \int_0^1 |\nabla u(x + t (y-x))|^p \, dt \,  \rho_n (|y-x|) \, d x \, d y.  
\end{multline}
By a change of variables $(x, z) = (x, x-y)$ and by applying the Fubini theorem, we obtain \eqref{thm-BBM-cl1} from \eqref{thm-BBM-p1}. We next deal with \eqref{thm-BBM-cl1} for which $\Phi(u) < + \infty$.  
The proof in this case follows from the previous case by considering a sequence $(u_k) \subset  C^\infty_c(\mR^N)$ such that $\Phi(u_k) \to \Phi(u)$ and $u_k \to u$ for almost every $x \in \mR^N$, as $k \to + \infty$. By Fatou's lemma, one has, for $n \ge 1$,  
\be \label{thm-BBM-p2}
\liminf_{k \to + \infty} \mathop{\iint}_{\mR^N \times \mR^N} \frac{|u_k(x)
- u_k(y)|^p}{|x-y|^p} \rho_n(|x-y|)\, dx \, dy 
\ge \mathop{\iint}_{\mR^N \times \mR^N} \frac{|u(x)
- u(y)|^p}{|x-y|^p} \rho_n(|x-y|)\, dx \, dy. 
\ee
Assertion \eqref{thm-BBM-cl1} now follows from \eqref{thm-BBM-p2} by applying \eqref{thm-BBM-p1} to $u_k$ and then letting $k \to + \infty$.

We now address \eqref{thm-BBM-cl2} for \hmr{$u$ such that} $\Phi(u)< + \infty$. Using the properties of the mollifier sequence $(\rho_n)$, in particular, the mass of $\rho_n$ concentrates around 0,  and \hmr{a} Taylor expansion, one can prove \eqref{thm-BBM-cl2} if $u \hmr{\in} C^\infty_c (\mR^N)$ in addition. Combining this with \eqref{thm-BBM-cl1}, one can derive \eqref{thm-BBM-cl2} for $u$ satisfying $\Phi(u) < + \infty$ after using the fact that, for $u_k, u \in L^p(\mR^N)$,  
\begin{multline*}
\left| \mathop{\iint}_{\mR^N \times \mR^N} \frac{|u_k(x)
- u_k(y)|^p}{|x-y|^p} \rho_n(|x-y|)\, dx \, dy - \mathop{\iint}_{\mR^N \times \mR^N} \frac{|u(x)
- u(y)|^p}{|x-y|^p} \rho_n(|x-y|)\, dx \, dy \right| \\[6pt]
\le C_p \mathop{\iint}_{\mR^N \times \mR^N} \frac{|(u_k - u)(x)
- (u_k-u) (y)|^p}{|x-y|^p} \rho_n(|x-y|)\, dx \, dy. 
\end{multline*}

It remains to prove that if the $\liminf_{n\to +\infty}$ of the LHS of \eqref{thm-BBM-cl1} is finite then $\Phi(u) < + \infty$. One way to prove this fact is to use the convexity of $t^p$ and the convolution suggested by Stein and presented in \cite{Brezis-02} as follows. Let $(\varphi_k)_{k \ge 1}$ be a smooth non-negative sequence of approximations to the identity such that $\supp \varphi_k \subset B_{1/k}$. Set 
$$
u_k = \varphi_k  * u \mbox{ for } k \ge 1. 
$$
Using Jensen's inequality, we derive that, for $n \ge 1$,  
$$
\mathop{\iint}_{\mR^N \times \mR^N} \frac{|u_k(x)
- u_k(y)|^p}{|x-y|^p} \rho_n(|x-y|)\, dx \, dy  \le \mathop{\iint}_{\mR^N \times \mR^N} \frac{|u(x)
- u(y)|^p}{|x-y|^p} \rho_n(|x-y|)\, dx \, dy \mbox{ for } k \ge 1. 
$$
Since $u_k \in C^\infty(\mR^N)$, one can show that, for $R>2$ and $0< r < 1$,  
\be
\lim_{n \to + \infty} \mathop{\iint}_{x  \in B_R; |y - x| < r} \frac{|u_k(x)
- u_k(y)|^p}{|x-y|^p} \rho_n(|x-y|)\, dx \, dy \ge K_{N, p} \int_{B_{R}} |\nabla u_k|^p \, dx. 
\ee
Here and in what follows,  $B_R$ denotes the open ball in $\mR^N$ centered at $0$ and of radius $R$ for $R>0$.  This implies 
$$
\lim_{n \to + \infty} \mathop{\iint}_{\mR^N \times \mR^N} \frac{|u_k(x)
- u_k(y)|^p}{|x-y|^p} \rho_n(|x-y|)\, dx \, dy \ge K_{N, p} \int_{\mR^N} |\nabla u_k|^p \, dx. 
$$
We thus obtain 
\be \label{thm-BBM-p3}
+ \infty > \liminf_{n \to + \infty} \mathop{\iint}_{\mR^N \times \mR^N} \frac{|u(x)
- u(y)|^p}{|x-y|^p} \rho_n(|x-y|)\, dx \, dy \ge  K_{N, p} \int_{\mR^N} |\nabla u_k|^p \, dx \mbox{ for } k \ge 1. 
\ee
Since $k \ge 1$ is arbitrary,  we obtain $\Phi(u) < + \infty$. The details of these arguments can be found in \cite{BBM-01} and \cite{Brezis-02}.

\medskip 
We next give a useful consequence of \Cref{thm-BBM}. 

\begin{proposition} \label{pro-BBM}
Let $N \ge 1$, $1 \le p < + \infty$, and $r > 0$, and let $u \in L^p(\mR^N)$. Then 
\begin{equation}\label{pro-BBM-cl1}
\lim_{\eps \to  0}  \mathop{\mathop{\iint}_{\mR^N \times \mR^N}}_{|x - y| < r} \frac{|u(x)
-u(y)|^p}{|x-y|^{p}} \frac{\eps}{|x-y|^{N-\eps}} \, dx \, dy =  K_{N, p} \Phi(u). 
\end{equation}
Consequently, if $u \in L^p(\mR^N)$ satisfies 
$$
\mathop{\iint}_{\mR^N \times \mR^N} \frac{|u(x)
-u(y)|^p}{|x-y|^{N+p}} \, dx \, dy <  + \infty, 
$$
then $u$ is constant. 
\end{proposition}

The reader can find many other interesting examples in \cite{BBM-01,Brezis-02,BSVY-222} on the way to determine whether or not a given function is constant. One of the motivations for determining whether or not a function is a constant is from the study of the Ginzburg-Landau equations, see, e.g., \cite{BBM-94}. Further properties related to the BBM formula can be found in \cite{BrNg-16,PS-17} (see also \cite{LS-11}) and the references therein. 

\begin{remark} \rm The limit of the term in the LHS of \eqref{pro-BBM-cl1} for $\eps \to 1_-$ was studied 
by Maz’ya and Shaposhnikova \cite{MS-02}. 
\end{remark}

\medskip 

We next discuss a related result of \Cref{thm-BBM} due to Nguyen \cite{NgSob1}, and  Bourgain and Nguyen \cite{BN-06}. 

\begin{theorem} [Bourgain \& Nguyen] \label{thm-I}
Let $N \ge 1$ and $1\le  p < + \infty$. The following two facts hold: 
\begin{enumerate}
\item[$i)$] For $p>1$ and  $u \in W^{1, p}(\mR^N)$, we have
\begin{equation} \label{thm-I-cl1-1}
\mathop{\mathop{\iint}_{\mR^N  \times \mR^N} }_{|u(x)-u(y)| > \delta}
\frac{\delta^p}{|x-y|^{N+p}} \, dx \, dy \le C_{N,p} \Phi(u)  \quad \mbox{ for } \delta > 0, 
\end{equation}
and 
\begin{equation} \label{thm-I-cl1-2}
\lim_{\delta \to 0_+} \mathop{\mathop{\iint}_{\mR^N  \times \mR^N} }_{|u(x)-u(y)| > \delta}
\frac{\delta^p}{|x-y|^{N+p}} \, dx \, dy =  \frac{1}{p}K_{N,p}  \Phi (u), 
\end{equation}
where $K_{N,p}$ is defined by \eqref{def-K}.

\item[$ii)$] If $u \in L^p(\mR^N)$ and $p \ge 1$, 
then 
\be \label{thm-I-cl2}
\Phi (u) \le C_{N, p}\liminf_{\delta \to 0_+}\mathop{\mathop{\iint}_{\mR^N  \times \mR^N} }_{|u(x)-u(y)| > \delta}
\frac{\delta^p}{|x-y|^{N+p}} \, dx \, dy. 
\ee
In particular, $u \in W^{1, p}(\mR^N)$ if $p>1$ and $u \in BV(\mR^N)$ if $p=1$ if the RHS of \eqref{thm-I-cl2} is finite. 
\end{enumerate}
\end{theorem}

\begin{remark} \rm The quantity given in the LHS of \eqref{thm-I-cl1-1} has its roots in estimates for the topological degree of continuous maps from a sphere into itself, which is due to Bourgain,  Brezis, and  Nguyen \cite{BBN-05} (see also \cite{BBM-04, BBM-05,Ng-Opt,Ng-refine}). 
\end{remark}

Assertion \eqref{thm-I-cl1-1} is based on the theory of maximal function. As in the analysis of \eqref{thm-BBM-cl1}, it suffices to consider the case $u \in C^\infty_c (\mR^N)$. Using spherical coordinates,  we have 
\be \label{thm-I-p0}
\mathop{\mathop{\iint}_{\mR^N  \times \mR^N} }_{|u(x)-u(y)| > \delta}
\frac{\delta^p}{|x-y|^{N+p}} \, dx \, dy = \mathop{\int_{\mR^N} \int_{\mS^{N-1}} \int_0^\infty }_{|u(x+ h \sigma) - u(x)| > \delta} \frac{\delta^p}{h^{p+1}} \, dh \, d \sigma \, dx. 
\ee
By a change of variables, we obtain 
\be  \label{thm-I-p0-111}
 \mathop{\int_{\mR^N} \int_{\mS^{N-1}} \int_0^\infty }_{|u(x+ h \sigma) - u(x)| > \delta} \frac{\delta^p}{h^{p+1}} \, dh \, d \sigma \, dx = \mathop{\int_{\mR^N} \int_{\mS^{N-1}} \int_0^\infty }_{\frac{1}{\delta h}|u(x+ \delta h \sigma) - u(x)| h > 1} \frac{1}{h^{p+1}} \, dh \, d \sigma \, dx,  
\ee
\hmr{Combining \eqref{thm-I-p0} and \eqref{thm-I-p0-111} yields 
\be \label{thm-I-p0-222}
\mathop{\mathop{\iint}_{\mR^N  \times \mR^N} }_{|u(x)-u(y)| > \delta}
\frac{\delta^p}{|x-y|^{N+p}} \, dx \, dy = \mathop{\int_{\mR^N} \int_{\mS^{N-1}} \int_0^\infty }_{\frac{1}{\delta h}|u(x+ \delta h \sigma) - u(x)| h > 1} \frac{1}{h^{p+1}} \, dh \, d \sigma \, dx. 
\ee}
Since  
$$
\frac{1}{\delta}|u(x+ \delta h \sigma) - u(x)| = h \left|\int_0^1 \nabla u (x + t \delta h \sigma) \cdot \sigma \, dt \right|  
\le M (\nabla u, \sigma) (x) h, 
$$
where 
$$
M(\nabla u, \sigma)(x): = \sup_{t > 0} \fint_{0}^t  \Big| \nabla u (x + s \sigma) \cdot \sigma| \, ds \hmr{,}
$$
\hmr{it} follows \hmr{\eqref{thm-I-p0-222}} that 
\begin{multline} \label{thm-I-p1-*}
\mathop{\mathop{\iint}_{\mR^N  \times \mR^N} }_{|u(x)-u(y)| > \delta}
\frac{\delta^p}{|x-y|^{N+p}} \, dx \, dy \le  \mathop{\int_{\mR^N} \int_{\mS^{N-1}} \int_0^\infty }_{M (\nabla u, \sigma) (x) h > 1} \frac{1}{h^{p+1}} \, dh \, d \sigma \, dx \\[6pt]= \frac{1}{p}  \int_{\mS^{N-1}} \int_{\mR^N} |M(\nabla u, \sigma)|^p \, dx \, d \sigma. 
\end{multline}
Applying the theory of maximal functions, see, e.g., \cite{Stein-70}, we have 
\be \label{thm-I-p2-*}
 \int_{\mR^N} |M(\nabla u, \sigma)|^p \, dx  =  \int_{\mR_\sigma^\perp}  \int_{\mR_\sigma}|M(\nabla u, \sigma)|^p \, dx 
 \le c_p \int_{\mR_\sigma^\perp} \int_{\mR_\sigma} |\nabla u \cdot \sigma|^p \, dx  \le c_p \int_{\mR^N} |\nabla u|^p. 
\ee
Here, for $\sigma \in \mS^{N-1}$,  
$$
\mR_\sigma = \{ t \sigma; t \in \mR \} \quad \mbox{ and } \quad \mR_\sigma^\perp = \{ x \in \mR^{N-1}; x \cdot \sigma = 0 \}. 
$$
Assertion \eqref{thm-I-cl1-1} for $u \in C^\infty_c(\mR^N)$ now follows from \eqref{thm-I-p1-*} and \eqref{thm-I-p2-*}. 

\medskip 
The proof of \eqref{thm-I-cl1-2} for $u \in C^\infty_c(\mR^N)$ follows from \hmr{\eqref{thm-I-p0-222}} and \hmr{a} Taylor expansion. The arguments used to prove \eqref{thm-I-cl1-2} in the general case \hmr{follow} from this case and \eqref{thm-I-cl1-1} as in the analysis of the BBM formula. 

\begin{remark} \rm Part $i)$ of \Cref{thm-I} was also obtained independently by Ponce and Van Schaftingen with a different proof as mentioned in \cite{NgSob1}.   
\end{remark}

We next briefly discuss the proof of \eqref{thm-I-cl2} under the stronger assumption in which the liminf is replaced by the limsup, i.e., 
\be \label{thm-I-cl2-2}
\Phi (u) \le C_{N, p}\limsup_{\delta \to 0}\mathop{\mathop{\iint}_{\mR^N  \times \mR^N} }_{|u(x)-u(y)| > \delta}
\frac{\delta^p}{|x-y|^{N+p}} \, dx \, dy. 
\ee
The ideas of the proof of \eqref{thm-I-cl2} under the assumption stated in \Cref{thm-I} will be later given in \Cref{sect-thm-neg-1p}. 

To be able to apply the arguments used in the proof of the BBM formula, one arranges to gain some convexity so that the arguments involving the convolution can be used. This can be done by an appropriate integration with respect to $
\delta $ as suggested \hmr{by the author} in \cite{NgSob1}. Indeed, we first assume that $u \in L^\infty(\mR^N)$ and let $a$ be a positive constant which is greater than $2 \| u\|_{L^\infty(\mR^N)} + 1$.  
We have, by using Fubini's theorem,  
\be \label{thm-I-p0-1}
\int_0^a \eps \delta^{-1 + \eps} \mathop{\mathop{\iint}_{\mR^N  \times \mR^N} }_{|u(x)-u(y)| > \delta}
\frac{\delta^p}{|x-y|^{N+p}} \, dx \, dy = \mathop{\iint}_{\mR^N  \times \mR^N} \frac{\eps}{p+\eps}
\frac{|u(x) - u(y)|^{p+\eps}}{|x-y|^{N+p}} \, dx \, dy.
\ee 
Let $(\varphi_k)_{k \ge 1}$ be a smooth non-negative sequence of approximations to the identity such that $\supp \varphi_k \subset B_{1/k}$. Set 
$$
u_k = \varphi_k  * u \mbox{ for } k \ge 1. 
$$
Since, by Jensen's inequality,  
$$
\mathop{\iint}_{\mR^N  \times \mR^N} \frac{\eps}{p+\eps}
\frac{|u_k(x) - u_k(y)|^{p+\eps}}{|x-y|^{N+p}} \, dx \, dy \le \mathop{\iint}_{\mR^N  \times \mR^N} \frac{\eps}{p+\eps}
\frac{|u(x) - u(y)|^{p+\eps}}{|x-y|^{N+p}} \, dx \, dy, 
$$
it follows from \eqref{thm-I-p0-1} that 
\be \label{thm-I-p1}
\int_0^a \eps \delta^{-1 + \eps} \mathop{\mathop{\iint}_{\mR^N  \times \mR^N} }_{|u(x)-u(y)| > \delta}
\frac{\delta^p}{|x-y|^{N+p}} \, dx \, dy \ge \mathop{\iint}_{\mR^N  \times \mR^N} \frac{\eps}{p+\eps}
\frac{|u_k(x) - u_k(y)|^{p+\eps}}{|x-y|^{N+p}} \, dx \, dy. 
\ee 

Using the fact $u_k \in C^\infty(\mR^N)$, we obtain, as in the proof of \eqref{thm-BBM-p3},  
$$
\mathop{\iint}_{\mR^N  \times \mR^N} \frac{\eps}{p+\eps}
\frac{|u_k(x) - u_k(y)|^{p+\eps}}{|x-y|^{N+p}} \, dx \, dy \ge C_{N, p} \int_{\mR^N} |\nabla u_k|^p \, dx, 
$$
which yields, by \eqref{thm-I-p1},  
$$
\limsup_{\delta \to 0_+} \mathop{\mathop{\iint}_{\mR^N  \times \mR^N} }_{|u(x)-u(y)| > \delta}
\frac{\delta^p}{|x-y|^{N+p}} \, dx \, dy  \ge C_{N, p} \int_{\mR^N} |\nabla u_k|^p \, dx. 
$$
Since $k \ge 1$ is arbitrary, the conclusion follows. 

The proof in the general case, without the assumption $u \in L^\infty(\mR^N)$, follows from this case after noting that 
\be \label{thm-I-p1-2}
\mathop{\mathop{\iint}_{\mR^N  \times \mR^N} }_{|u(x)-u(y)| > \delta}
\frac{\delta^p}{|x-y|^{N+p}} \, dx \, dy \ge \mathop{\mathop{\iint}_{\mR^N  \times \mR^N} }_{|u_A(x)-u_A(y)| > \delta}
\frac{\delta^p}{|x-y|^{N+p}} \, dx \, dy,  
\ee
where 
$$
u_A = \min\{ \max\{u, - A\}, A\}. 
$$
Indeed, since $u_A \in L^\infty(\mR^N)$, we have 
\be \label{thm-I-p2}
\limsup_{\delta \to 0_+ }\mathop{\mathop{\iint}_{\mR^N  \times \mR^N} }_{|u_A(x)-u_A(y)| > \delta}
\frac{\delta^p}{|x-y|^{N+p}} \, dx \, dy \ge C_{N, p} \Phi(u_A). 
\ee
Combining \eqref{thm-I-p1-2} and \eqref{thm-I-p2} yields 
$$
\mathop{\mathop{\iint}_{\mR^N  \times \mR^N} }_{|u(x)-u(y)| > \delta}
\frac{\delta^p}{|x-y|^{N+p}} \, dx \, dy \ge C_{N, p} \Phi(u_A), 
$$
which implies \eqref{thm-I-cl2-2} after letting $A \to + \infty$.  

\medskip 
The setting considered in \Cref{thm-I} has been extended in several directions. See, e.g., \cite{NgSob2,BN-18,BVY-21,BSVY-24} where a more general functional was considered, see, e.g., \cite{NPSV-18,PSV-19} where a magnetic field was involved  (related properties for the BBM formula \hmr{were} also considered, see, e.g., \cite{NPSV-18,PSV-19,NS-19}). 

\medskip
One of these extensions, recently proposed by Brezis, Seeger, Van Schaftingen, and Yung \cite{BSVY-24}, will be now discussed.  
Set, for a measurable function $u$ defined in $\mR^N$,  
$$
Q_b u (x, y) = \frac{u(x) - u(y)}{|x - y|^{1 + b}} \mbox{ for } x, y \in \mR^N, 
$$
and, for a measurable set $E$ of $\mR^{N} \times \mR^{N}$, 
$$
\nu_{\gamma} (E) = \iint_{E} |x - y|^{-N+ \gamma} \, dx \, dy, 
$$
and 
$$
 E_{\lambda, b} (u)= \Big\{ (x, y); |Q_b u(x, y)| > \lambda \Big\}. 
$$

Given $\gamma \in \mR$ and $p \ge 1$, set 
$$
\Phi_\lambda(u) = \nu_{\gamma} \Big( E_{\lambda, \gamma/p} (u) \Big) \mbox{ for } \lambda > 0. 
$$
Given a measurable subset $\Omega$ of $\mR^N$ and a measurable function $u$ defined in $\Omega$, we also denote  
$$
 E_{\lambda, b, \Omega} (u)= \Big\{ (x, y) \in \Omega \times \Omega; |Q_b u(x, y)| > \lambda \Big\}. 
$$
and 
$$
\Phi_{\lambda, \Omega} (u) = \lambda^p \nu_{\gamma} \Big( E_{\lambda, \gamma/p, \Omega} (u) \Big). 
$$

One has \hmr{the following result.}

\begin{theorem}[Brezis \& Seeger \& Van Schaftingen \& Yung] \label{thm-B} Let $N \ge 1$, $1 \le   p < + \infty$, $\gamma \in \mR \setminus \{0 \}$. The following facts hold: 
\begin{enumerate}
\item[$i)$] For $p>1$ and  $u \in W^{1, p}(\mR^N)$, we have
\begin{equation} \label{thm-B-cl1-1}
\Phi_\lambda (u) \le C_{N,p}  \Phi(u)  \quad  \mbox{ for } \lambda > 0, 
\end{equation}
and 
\begin{equation} \label{thm-B-cl1-2}
\lim_{\lambda \to 0_+} \Phi_\lambda (u) =  \frac{1}{|\gamma|}K_{N,p} \Phi (u) \mbox{ if } \gamma < 0, 
\end{equation}
and 
\begin{equation} \label{thm-B-cl1-3}
\lim_{\lambda \to + \infty} \Phi_\lambda (u) =  \frac{1}{|\gamma|}K_{N,p} \Phi(u) \mbox{ if } \gamma > 0, 
\end{equation}
where $K_{N,p}$ is defined by \eqref{def-K}.

\item[$ii)$] If $u \in L^p(\mR^N)$ and $p \ge 1$, then 
\be \label{thm-B-cl2}
\Phi(u) \le C \sup_{\lambda > 0} \Phi_\lambda (u). 
\ee
In particular, $u \in W^{1, p}(\mR^N)$ if $p>1$ and $u \in BV(\mR^N)$ if $p=1$ if the RHS of \eqref{thm-B-cl2} is finite. 
\end{enumerate}

\end{theorem}

\begin{remark} \rm For $\gamma = -p$, one has 
$$
\Phi_\delta (u) = \mathop{\mathop{\iint}_{\mR^N  \times \mR^N} }_{|u(x)-u(y)| > \delta}
\frac{\delta^p}{|x-y|^{N+p}} \, dx \, dy, 
$$
the quantity considered in \Cref{thm-I}. For $\gamma = N$, the quantity $\Phi_\lambda$ was previously considered by Brezis, Van Schaftingen, and Yung \cite{BVY-21}. 
A generalization in this case to
one-parameter families of operators was considered by Domínguez and Milman \cite{DM-22}.
\end{remark}

\begin{remark} \rm It is worth noting that \eqref{thm-B-cl1-2} and \eqref{thm-B-cl1-3} fail for $p=1$ and for $u \in BV(\mR^N)$ or  even for $u \in W^{1, 1}(\mR^N)$. This has been noted by Brezis, Seeger, Van Schaftingen, and Yung  \cite{BSVY-24}. This phenomenon in the case $\gamma = -1$ was previously observed by Ponce, see also the work of Brezis and Nguyen \cite{BN-18}. 
\end{remark}

We next briefly discuss the proof of part $i)$ of \Cref{thm-B}.  The proof is closely related to the proof of \Cref{thm-I} mentioned previously and is different from the one given in \cite{BSVY-24}. We only discuss the case $N=1$. The general case \hmr{can be established} as in the proof of \Cref{thm-I}.  We begin with the proof of assertion \eqref{thm-B-cl1-1}. As in the proof of \eqref{thm-I}, it suffices to consider only the case $u \in C^\infty_c(\mR)$.  
We have 
\be \label{thm-B-pp0}
 \mathop{\mathop{\iint}_{\mR \times \mR}}_{\frac{|u(x) - u(y)|}{|x-y|^{1 + \frac{\gamma}{p}}} > \lambda} |x-y|^{-1 + \gamma} \, dx \, dy = 2  \mathop{\int_{\mR} \int_0^\infty}_{\frac{|u(x+ h) - u(h)|}{h^{1 + \frac{\gamma}{p}}} > \lambda} h^{-1 + \gamma} \, dh  \, dx. 
\ee
Since 
$$
|u(x+ h) - u(x)| = \left| \int_{x}^{x+ h} u'(s) \, ds \right| \le  \int_{x}^{x+ h} |u'(s)| \, ds \le M(u')(x) h,  
$$
where $M(u')(x) = \sup_{h > 0} \fint_{x}^{x+h} |u'(s)| \, ds$, it follows that 
\be
 \mathop{\int_{\mR} \int_0^\infty}_{\frac{|u(x+ h) - u(x)|}{h^{1 + \frac{\gamma}{p}}} > \lambda} h^{-1 + \gamma} \, dh  \, dx \le  \mathop{\int_{\mR} \int_0^\infty}_{\frac{M(u')(x)}{h^{\frac{\gamma}{p}}} > \lambda} h^{-1 + \gamma} \, dh  \, dx.
\ee
By considering $\gamma > 0$ and $\gamma < 0$ separately, one can show that 
\be
\lambda^p \mathop{\int_{\mR} \int_0^\infty}_{\frac{M(u')(x)}{h^{\frac{\gamma}{p}}} > \lambda} h^{-1 + \gamma} \, dh  \, dx \le C_\gamma \int_{\mR}  |M(u')(x)|^p \, dx. 
\ee
Since, by the theory of maximal functions, 
$$
\int_{\mR} |M(u')(x)|^p \, dx \le c_p \int_{\mR} |u'|^p \, dx, 
$$
assertion \eqref{thm-B-cl1-1} follows for $u \in C^\infty_c(\mR)$. 

To prove \eqref{thm-B-cl1-2} and \eqref{thm-B-cl1-3}, one notes that, by a change of variables, with $\delta^{- \frac{\gamma}{p}} = \lambda$, 
$$
\lambda^p \mathop{\int_{\mR} \int_0^\infty}_{\frac{|u(x+ h) - u(x)|}{h^{1 + \frac{\gamma}{p}}} > \lambda} h^{-1 + \gamma} \, dh  \, dx =  \mathop{\int_{\mR} \int_0^\infty}_{\frac{|u(x+ \delta h) - u(x)|}{\delta h} h^{-\frac{\gamma}{p}} > 1} h^{-1 + \gamma} \, dh  \, dx,  
$$
which yields, by \eqref{thm-B-pp0}, with $\delta^{- \frac{\gamma}{p}} = \lambda$,  
\be \label{thm-B-p1}
\Phi_\lambda(u) =  \mathop{\int_{\mR} \int_0^\infty}_{\frac{|u(x+ \delta h) - u(x)|}{\delta h} h^{-\frac{\gamma}{p}} > 1} h^{-1 + \gamma} \, dh  \, dx. 
\ee

The proof of \eqref{thm-B-cl1-2} and  \eqref{thm-B-cl1-3} for $u \in C^\infty_c(\mR)$ follows from \eqref{thm-B-p1} and \hmr{a} Taylor expansion. The arguments in the general case follows from this case and \eqref{thm-B-cl1-1} as in the analysis of the BBM formula. 

The proof of part $ii)$ of \Cref{thm-B} given in \cite{BSVY-24} is based on the BBM formula and the Lorentz duality, which in turn involves 
rearrangement properties.  Later, we state and prove a stronger version of part $ii)$ (see \Cref{sect-Sobolev}). Our proof is in the spirit of the one of \Cref{thm-I} and thus different from \cite{BSVY-24}.

\medskip 
Brezis, Seeger, Van Schaftinger, and Yung \cite{BSVY-24} also proved that, for $p=1$,  the following result. 

\begin{proposition}[Brezis \& Seeger \& Van Schaftingen \& Yung]  \label{pro-B} Let $N \ge 1$, $p=1$, and $\gamma \in \mR \setminus \{0 \}$. Then 
\eqref{thm-B-cl1-2} holds for $u \in C^1_c(\mR^{\hmr{N}})$ if $\gamma \not \in [-1, 0)$ and for $u \in W^{1, 1}(\mR^N)$ if $\gamma > 0$. We also have, for $\gamma \not \in [-1, 0]$ 
\begin{equation} \label{pro-B-cl1}
\Phi_\lambda (u) \le C_{N,1} \| \nabla u \|_{\cM},  \quad \forall \, \lambda > 0, 
\end{equation}
\end{proposition}

\begin{remark} \label{rem-gamma-0-1} \rm Brezis \& Seeger \& Van Schaftingen \& Yung \cite{BSVY-24} also showed that \eqref{pro-B-cl1} fails for $ \gamma \in [-1, 0)$. 
\end{remark}

\begin{remark} \rm Assertion \eqref{thm-B-cl1-1} 
for $\gamma = N$ was obtained by Brezis, Van Schaftinger, and Yung \cite{BVY-21b}. They also discussed there variants of  Gagliardo \& Nirenberg interpolation inequalities for functions in $W^{1,1}(\mR^N)$.  
\end{remark}

The proof of \eqref{pro-B-cl1} is based on the Vitali covering lemma in the case $\gamma > 0$ and involves a clever way to estimate double integrals in the case $\gamma < -1$ (and $N=1$), see \cite{BSVY-24}.

\medskip 
Viewing \Cref{thm-B} and \Cref{pro-B}, the following questions are proposed, in the spirit of \Cref{thm-I}. 

\begin{question}\label{question-limsup} Let $N \ge 1$, $p \ge 1$, and  $\gamma \in \mR \setminus \{0\}$, and let $u \in L^p(\mR^N)$. Is it true that, for some positive constant $C_{N, p}$,  
$$
\limsup_{\lambda \to 0_+} \Phi_\lambda (u) \ge C_{N, p} \Phi (u) \mbox{ for } \gamma < 0,
$$
and
$$
\limsup_{\lambda \to + \infty} \Phi_\lambda (u) \ge C_{N, p} \Phi (u) \mbox{ for } \gamma > 0. 
$$
\end{question}

A stronger version which have been proposed by Brezis \& Seeger \& Van Schaftingen \& Yung \cite{BSVY-24} is the following question. 

\begin{question} \label{question-liminf} Let $N \ge 1$, $p \ge 1$, and  $\gamma \in \mR \setminus \{0\}$, and let $u \in L^p(\mR^N)$. Is it true that, for some positive constant $C_{N, p}$,  
$$
\liminf_{\lambda \to 0_+} \Phi_\lambda (u) \ge C_{N, p} \Phi (u) \mbox{ for } \gamma < 0
$$
and
$$
\liminf_{\lambda \to + \infty} \Phi_\lambda (u) \ge C_{N, p} \Phi (u) \mbox{ for } \gamma > 0. 
$$
\end{question}

If the answer to  \Cref{question-liminf}  or \Cref{question-limsup} is positive, one then improves  \eqref{thm-B-cl1-1} and \eqref{pro-B-cl1}. 
Various results related to these questions will be stated in the next sections, \Cref{sect-Sobolev} - \Cref{sect-thm-neg-1p}. In \Cref{sect-Gamma}, we discuss about the Gamma-convergence of $\Phi_\lambda$,  and in \Cref{sect-Inequality} we discuss various inequalities related to $\Phi_\lambda$. 

\section{Characterizations of the Sobolev norms and the total variations} \label{sect-Sobolev}

In this section, we present various positive and negative results related to \Cref{question-limsup} and \Cref{question-liminf}. These results \hmr{particularly} \hmr{improve} part $ii)$ of \Cref{thm-B}. Here are the main results in this direction. We begin with the case $\gamma < 0$.

\begin{theorem} \label{thm-neg-p} Let $N \ge 1$, $1 \le p  < + \infty$, and $\gamma < -p$, and  let $u \in L^p(\mR^n)$. We have 
\be \label{thm-neg-p-cl}
\frac{K_{N, p}}{|\gamma|} \Phi(u)  \le \limsup_{\lambda \to 0_+} \Phi_\lambda (u).  
\ee
In particular, $u \in W^{1, p}(\mR^N)$ if $p> 1$ and $u \in BV(\mR^N)$ if $p=1$ if the RHS of \eqref{thm-neg-p-cl} is finite. 
\end{theorem} 

The proof of \Cref{thm-neg-p} is in the spirit of the one of \eqref{thm-I-cl2-2} and is given in \Cref{sect-thm-neg-p}.

\begin{theorem} \label{thm-neg-1p} Let $N \ge 1$, $1 \le p  < + \infty$, and $-p \le \gamma \le - 1$, and let $u \in L^p(\mR^N)$. We have
\be \label{thm-neg-1p-cl}
C_{N, p} \Phi(u) \le  \liminf_{\lambda \to 0_+} \Phi_\lambda (u)
\ee
for some positive constant $C_{N, p}$. In particular, $u \in W^{1, p}(\mR^N)$ if $p> 1$ and $u \in BV(\mR^N)$ if $p=1$ if the RHS of \eqref{thm-neg-1p-cl} is finite. 
\end{theorem}

As mentioned previously, the case $\gamma = - p$ is due to Bourgain and Nguyen \cite{BN-06}. The proof of \Cref{thm-neg-p} is in the same spirit and is presented in \Cref{sect-thm-neg-1p}.

\medskip 
The case $-1 < \gamma < 0$ is quite special. We have the following result whose proof  is given in \Cref{sect-pro-CE}. 

\begin{proposition}\label{pro-CE} Let $N \ge 1$, $1 \le p  < + \infty$, and $-1 < \gamma < 0$. There exists \hmr{a non-zero function} $u \in BV(\mR^N)$ with compact support such that 
$$
\lim_{\lambda \to 0_+} \Phi_\lambda (u) = 0. 
$$
\end{proposition} 

\Cref{pro-CE} gives a negative answer to \Cref{question-limsup} and \Cref{question-liminf} in the case $-1 < \gamma < 0$. 
Nevertheless, we can prove the following result which improves part $ii)$ of \Cref{thm-B} in the case $-1 < \gamma < 0$.

\begin{theorem}\label{thm-pos-neg} Let $N \ge 1$, $1 \le p  < + \infty$, and $-1 < \gamma < 0$, and  let $u \in L^p(\mR^n)$.  Then 
\be \label{thm-pos-neg-cl}
 \frac{K_{N, p}}{\gamma + \hmr{2} p } \Phi(u) \le  \hmr{\limsup_{\lambda \to 0_+} \Phi_\lambda (u)} +  \limsup_{\lambda \to + \infty} \Phi_\lambda (u). 
\ee 
In particular, $u \in W^{1, p}(\mR^N)$ if $p> 1$ and $u \in BV(\mR^N)$ if $p=1$ if the RHS of \eqref{thm-pos-neg-cl} is finite. 
\end{theorem}

\hmr{The proof of \Cref{thm-pos-neg} is in the spirit of the one of \Cref{thm-I} mentioned in the introduction and given in \Cref{sect-thm-pos-neg}.}

\medskip 

We next deal with the case $\gamma > 0$. 
\begin{theorem}\label{thm-pos} Let $N \ge 1$, $1 \le p  < + \infty$, and $\gamma > 0$, and  let $u \in L^1_{\loc}(\mR^n)$.  Then 
\be \label{thm-pos-cl}
 \frac{K_{N, p}}{\gamma + p } \Phi(u) \le  \limsup_{\lambda \to + \infty} \Phi_\lambda (u). 
\ee 
In particular, $u \in W^{1, p}(\mR^N)$ if $p> 1$ and $u \in BV(\mR^N)$ if $p=1$ if $u \in L^p(\mR^N)$ and the RHS of \eqref{thm-pos-cl} is finite. 
\end{theorem} 

The proof of \Cref{thm-pos}  is in the spirit of the one of \Cref{thm-I} mentioned in the introduction and  given  in \Cref{sect-thm-pos}. 

\medskip 
\Cref{thm-pos} is known in the case $p>1$ and partially known in the case $p=1$. In the case $\gamma = \hmr{N}$ and $p \ge 1$, \Cref{thm-pos} was established by Poliakovsky \cite{Poliakovsky-22} with $C_{N, p}$ instead of $\frac{K_{N, p}}{\gamma + p}$ in the LHS of \eqref{thm-pos-cl}. A stronger result of \Cref{thm-pos} with $C_{N, p, \gamma}$, \hmr{blowing up} as $\gamma \to 0_+$, instead of $\frac{K_{N, p}}{\gamma + p}$ in the LHS of \eqref{thm-pos-cl} but with liminf instead of limsup on the RHS was obtained by Gobbino and Picenni \cite{GB-25}. Their results do not imply ours in the case $p=1$ since \eqref{thm-B-cl1-3} does not hold for $p=1$. \Cref{thm-pos} for $p=1$ is sharp in the sense that the equality holds in \eqref{thm-pos-cl} when $u$ is the characteristic function of a bounded convex domain with smooth boundary, see \cite[Lemma 3.6]{BSVY-24} (see also \cite{Picenni-24}).

\section{Proof of \Cref{thm-neg-p}} \label{sect-thm-neg-p}
 The proof is based on the ideas of the convexity and the BBM formula as mentioned in the proof of \eqref{thm-I-cl2-2}. We first assume that $u \in L^\infty(\mR^N)$ and let $m \ge 1$ be such that 
\be
\| u\|_{L^\infty(\mR^N)} \le m. 
\ee

Fix $a > 0$ sufficiently small. Then 
\be \label{thm-neg-p-p0}
\limsup_{\lambda \to 0_+} \Phi_\lambda (u) \ge \limsup_{\eps \to 0_+} \int_0^a \frac{\eps}{\lambda^{1 - \eps}} \Phi_\lambda(u) \, d \lambda. 
\ee
We have, for $\eps > 0$,  
\begin{multline*}
\int_0^a \frac{\eps}{\lambda^{1 - \eps}} \Phi_{\lambda} (u) \, d \lambda = \int_0^a \frac{\eps}{\lambda^{1 - \eps}}  \lambda^p \mathop{\iint}_{\frac{|u(x) - u(y)|}{|x-y|^{1 + \frac{\gamma}{p}} } > \lambda} |x-y|^{-N + \gamma} \, dx \, dy \\[6pt]
= \mathop{\iint}_{\mR^N \times \mR^N}  |x-y|^{-N + \gamma} d x \, dy \int_0^a \eps \lambda^{p- 1 + \eps} \mathds{1}_{\lambda < \frac{|u(x) - u(y)|}{|x-y|^{1 + \frac{\gamma}{p}} }} \, d \lambda \\[6pt]
\ge  \mathop{\iint}_{\frac{|u(x) - u(y)|}{|x-y|^{1 + \frac{\gamma}{p}}} < a}|x-y|^{-N + \gamma}   \frac{\eps}{p+\eps}  \frac{|u(x) - u(y)|^{p+ \eps}}{|x-y|^{p+ \eps + \frac{\gamma (p+ \eps)}{p}}} d x \, dy. 
\end{multline*}
We thus obtain 
\be\label{thm-neg-p-p1}
\int_0^a \frac{\eps}{\lambda^{1 - \eps}} \Phi_{\lambda} (u) \, d \lambda 
\ge  \mathop{\iint}_{\frac{|u(x) - u(y)|}{|x-y|^{1 + \frac{\gamma}{p}}} < a}\frac{1}{p+\eps}   \frac{|u(x) - u(y)|^{p+ \eps}}{|x-y|^{p+ \eps}} \frac{\eps}{|x-y|^{N + \frac{\gamma \eps}{p}}}d x \, dy.
\ee

Set 
$$
\beta = 1 + \frac{\gamma}{p}. 
$$
Since $\gamma < -p$, it follows that $ \beta < 0.$  Denote
$$
r = (2m/a)^{1/\beta}.   
$$
Since $\beta < 0$, it follows that if $|x-y| < r$ then $\frac{|u(x) - u(y)|}{|x- y|^\beta} < \frac{2m}{r^\beta} = a$. 
We derive from \eqref{thm-neg-p-p1} that 
\be \label{thm-neg-p-p2}
\int_0^a \frac{\eps}{\lambda^{1 - \eps}} \Phi_{\gamma} (u) \, d \lambda 
\ge  \mathop{\iint}_{|x-y| < r} \frac{1}{p+\eps}   \frac{|u(x) - u(y)|^{p+ \eps}}{|x-y|^{p+ \eps}} \frac{\eps}{|x-y|^{N + \frac{\gamma \eps}{p}}}d x \, dy.
\ee

Let $(\varphi_k)$ be a smooth non-negative sequence of approximations to the identity such that $\supp \varphi_k \subset B_{1/k}$. Set 
$$
u_k = \varphi_k  * u \mbox{ for } k \ge 1. 
$$
Since $p + \eps >  1$, we have, by Jensen's inequality,
\begin{multline} \label{thm-neg-p-p3}
\mathop{\iint}_{|x-y| < r} \frac{1}{p+\eps}   \frac{|u(x) - u(y)|^{p+ \eps}}{|x-y|^{p+ \eps}} \frac{\eps}{|x-y|^{N + \frac{\gamma \eps}{p}}}d x \, dy \\[6pt] \ge \mathop{\iint}_{|x-y| < r} \frac{1}{p+\eps}   \frac{|u_k(x) - u_k(y)|^{p+ \eps}}{|x-y|^{p+ \eps}} \frac{\eps}{|x-y|^{N + \frac{\gamma \eps}{p}}}d x \, dy. 
\end{multline}

Combining  \eqref{thm-neg-p-p2} and  \eqref{thm-neg-p-p3} yields, for $k \ge 1$,  
\be
\int_0^a \frac{\eps}{\lambda^{1 - \eps}}  \Phi_{\gamma} (u) \, d \lambda \ge \mathop{\iint}_{|x-y| < r} \frac{1}{p+\eps}   \frac{|u_k(x) - u_k(y)|^{p+ \eps}}{|x-y|^{p+ \eps}} \frac{\eps}{|x-y|^{N + \frac{\gamma \eps}{p}}}d x \, dy. 
\ee
By letting $\eps \to 0_+$, and noting that $\gamma < 0$ and using \eqref{thm-neg-p-p0}, we obtain, as in the proof of \eqref{thm-BBM-p3},  
\be
\limsup_{\lambda \to 0} \Phi_\lambda (u) \ge \frac{1}{|\gamma|} K_{N,p} \int_{\mR^N} |\nabla u_k|^p \, dx
\ee
\hmr{since, for $\gamma < 0$, 
$$
\lim_{\eps \to 0_+} \int_0^r \frac{\eps}{s^{1 + \frac{\gamma \eps}{p}}} \, ds = \frac{p}{|\gamma|}. 
$$
}
By taking $k \to + \infty$, we reach the conclusion when $u \in L^\infty(\mR^N)$. The proof in the general case follows from the case where $u \in L^\infty(\mR^N)$ as in the proof of \Cref{thm-I} after noting that 
$$
\Phi_\lambda (u) \ge \Phi_\lambda (u_A), 
$$
where 
$$
u_A = \min\{ \max\{u, - A\}, A\}. 
$$
The details are omitted. \qed

\section{Proof of \Cref{thm-pos}} \label{sect-thm-pos}

The proof is again based on the ideas of the convexity and the BBM formula as mentioned in the proof of \eqref{thm-I-cl2-2} related to \Cref{thm-I}. We first assume that 
$$
\| u\|_{L^\infty} \le m. 
$$

Let $a > 0$ be sufficiently large. Then 
\be \label{thm-pos-neg-p0}
\limsup_{\lambda \to + \infty} \Phi_\lambda (u) \ge \limsup_{\eps \to 0_+} \int_a^\infty \frac{\eps}{\lambda^{1 + \eps}} \Phi_{\lambda} (u) \, d \lambda. 
\ee
We have 
\begin{multline*}
\int_a^\infty \frac{\eps}{\lambda^{1 + \eps}} \Phi_{\lambda} (u) \, d \lambda = \int_a^\infty \frac{\eps}{\lambda^{1 + \eps}}  \lambda^p \mathop{\iint}_{\frac{|u(x) - u(y)|}{|x-y|^{1 + \frac{\gamma}{p}} } > \lambda} |x-y|^{-N + \gamma} \, dx \, dy 
\\[6pt]
= \mathop{\iint}_{\mR^N \times \mR^N}  |x-y|^{-N + \gamma} d x \, dy \int_a^\infty \eps \lambda^{p- 1 - \eps} \mathds{1}_{\lambda < \frac{|u(x) - u(y)|}{|x-y|^{1 + \frac{\gamma}{p}} }} \, d \lambda. 
\end{multline*}
It follows that, for every $R> 2$ and $ 0< r < 1$, 
\begin{multline} \label{thm-pos-neg-p1}
\int_a^\infty \frac{\eps}{\lambda^{1 + \eps}} \Phi_{\lambda} (u) \, d \lambda \ge \mathop{\iint}_{x \in B_R; |x-y| < r}  |x-y|^{-N + \gamma} d x \, dy \int_a^\infty \eps \lambda^{p- 1 - \eps} \mathds{1}_{\lambda < \frac{|u(x) - u(y)|}{|x-y|^{1 + \frac{\gamma}{p}} }} \, d \lambda 
\\[6pt]
\ge  \mathop{\iint}_{x \in B_R; |x-y| < r}   |x-y|^{-N + \gamma}  \left( \frac{\eps}{p-\eps} \frac{|u(x) - u(y)|^{p- \eps}}{|x-y|^{p- \eps + \frac{\gamma (p- \eps)}{p}}} - \frac{\eps}{p-\eps} a^{p-\eps}  \right) d x \, dy 
\\[6pt]
= \frac{1}{p-\eps}\mathop{\iint}_{x \in B_R; |x-y| < r}     \frac{|u(x) - u(y)|^{p- \eps}}{|x-y|^{p- \eps}} \frac{\eps}{|x-y|^{N - \frac{\eps \gamma}{p}}}  d x \, dy \\[6pt]
 -  \frac{\eps a^{p-\eps}}{p-\eps}  \mathop{\iint}_{x \in B_R; |x-y| < r}  |x-y|^{-N + \gamma} \, dx \, dy. 
\end{multline}
We have 
\begin{multline}\label{thm-pos-neg-p2}
\frac{1}{p-\eps}\mathop{\iint}_{x \in B_R; |x-y| < r}    \frac{|u(x) - u(y)|^{p- \eps}}{|x-y|^{p -  \eps }}  \frac{\eps}{|x-y|^{N - \frac{\eps \gamma}{p}}}d x \, dy \\[6pt]
\ge \frac{1}{p-\eps} \frac{1}{|2m|^{\eps}}\mathop{\iint}_{x \in B_R\hmr{;}   |x-y| < r}     \frac{|u(x) - u(y)|^{p}}{|x-y|^{p}}  \frac{\eps}{|x-y|^{N - \frac{\eps \gamma}{p} - \eps}}d x \, dy. 
\end{multline}

Let $(\varphi_k)_{k \ge 1}$ be a smooth non-negative sequence of approximations to the identity such that $\supp \varphi_k \subset B_{1/k}$. Set 
$$
u_k = \varphi_k  * u \mbox{ for } k \ge 1. 
$$
We have, by Jensen's inequality,  
\begin{multline}\label{thm-pos-neg-p3}
 \frac{1}{p-\eps} \frac{1}{|2m|^{\eps}}\mathop{\iint}_{x \in B_R; |x-y| < r}    \frac{|u(x) - u(y)|^{p}}{|x-y|^{p - \eps}} \frac{\eps}{|x-y|^{N - \frac{\eps \gamma}{p} - \eps}} d x \, dy \\[6pt] \ge  \frac{1}{p-\eps} \frac{1}{|2m|^{\eps}} \mathop{\iint}_{x \in B_{R-1}; |x-y| < r}  \frac{|u_k(x) - u_k(y)|^{p}}{|x-y|^{p}} \frac{\eps}{|x-y|^{N - \frac{\eps \gamma}{p} - \eps}} d x \, dy. 
\end{multline}
Since, thanks to fact that $\gamma > \hmr{0}$, 
\be
\lim_{\eps \to 0_+} \frac{\eps a^{p-\eps}}{p-\eps}  \mathop{\iint}_{x \in B_R; |x-y| < r}  |x-y|^{-N + \gamma} \, dx \, dy = 0, 
\ee
we derive from \eqref{thm-pos-neg-p1}, \eqref{thm-pos-neg-p2}, and \eqref{thm-pos-neg-p3} that 
$$
\limsup_{\eps \to 0_+} \int_a^\infty \frac{\eps}{\lambda^{1 - \eps}} \Phi_{\lambda} (u) d \lambda \ge \liminf_{\eps \to 0_+} \frac{1}{p-\eps} \frac{1}{|2m|^{\eps}}\mathop{\iint}_{x \in B_{R-1}; |x-y| < 1}    \frac{|u_k(x) - u_k(y)|^{p}}{|x-y|^{p}} \frac{\eps}{|x-y|^{N - \frac{\eps \gamma}{p} - \eps}}  d x \, dy. 
$$
Using \eqref{thm-pos-neg-p0}, we obtain  
\be
\limsup_{\lambda \to + \infty} \Phi_\lambda(u)  \ge \frac{K_{N,p}}{\gamma + p}  \int_{B_{R-1}} |\nabla u_k|^p \, dx 
\ee
\hmr{since, for $-1 < \gamma < 0$, 
$$
\lim_{\eps \to 0_+} \int_0^1 \frac{\eps}{s^{1 - \frac{\gamma \eps}{p}  - \eps}} \, ds = \frac{p}{\gamma + p}. 
$$
}
Since $R>2$ is arbitrary, we reach 
\be
\limsup_{\lambda \to + \infty} \Phi_\lambda(u) \ge \frac{K_{N,p}}{\gamma + p}  \int_{\mR^N} |\nabla u_k|^p \, dx. 
\ee
Letting $k \to + \infty$, we reach the conclusion when $u \in L^\infty(\mR^N)$. 

The proof in the general case follows from the case where $u \in L^\infty(\mR^N)$ as in the proof of \Cref{thm-I} after noting that 
$$
\Phi_\lambda (u) \ge \Phi_\lambda (u_A), 
$$
where 
$$
u_A = \min\{ \max\{u, - A\}, A\}. 
$$
The details are omitted.
\qed

\section{Proof of \Cref{thm-pos-neg}} \label{sect-thm-pos-neg}

{\color{black}
The proof is again based on the ideas of the convexity and the BBM formula as mentioned in the proof of \eqref{thm-I-cl2-2} related to \Cref{thm-I}. Without loss of generality, one might assume that 
\be \label{thm-neg-p0-1-***}
\limsup_{\lambda \to 0_+} \Phi_\lambda u < + \infty. 
\ee
This implies, for all $1< \Lambda < + \infty$,  
\be \label{thm-neg-p0-1}
 \Phi_\lambda u < C_{\Lambda} \mbox{ for } \Lambda^{-1} \le \lambda \le \Lambda.  
\ee

We first assume that 
$$
\| u\|_{L^\infty} \le m \quad \mbox{ and } \quad m \ge 1. 
$$

Then, by \eqref{thm-neg-p0-1}, it holds 
\be \label{thm-neg-p0}
\limsup_{\lambda \to + \infty} \Phi_\lambda (u) \ge \limsup_{\eps \to 0_+} \int_1^\infty \frac{\eps}{\lambda^{1 + \eps}} \Phi_{\lambda} (u) \, d \lambda. 
\ee
We have, for every $R> 2$ and $ 0< r < 1$, 
\begin{multline} \label{thm-neg-p1}
\int_1^\infty \frac{\eps}{\lambda^{1 + \eps}} \Phi_{\lambda} (u) \, d \lambda 
\ge \frac{1}{p-\eps} \mathop{\mathop{\iint}_{x \in B_R; |x-y| < r}}_{\frac{|u(x) - u(y)|}{|x - y|^{1+ \frac{\gamma}{p}}} > 1}     \frac{|u(x) - u(y)|^{p- \eps}}{|x-y|^{p- \eps}} \frac{\eps}{|x-y|^{N - \frac{\eps \gamma}{p}}}  d x \, dy \\[6pt]
 -  \frac{\eps }{p-\eps}   \mathop{\mathop{\iint}_{x \in B_R; |x-y| < r}}_{\frac{|u(x) - u(y)|}{|x - y|^{1+ \frac{\gamma}{p}}} > 1}  |x-y|^{-N + \gamma} \, dx \, dy. 
\end{multline}
On the other hand,  
\begin{multline}\label{thm-neg-p2}
\frac{1}{p-\eps} \mathop{\mathop{\iint}_{x \in B_R; |x-y| < r}}_{\frac{|u(x) - u(y)|}{|x - y|^{1+ \frac{\gamma}{p}}} > 1}     \frac{|u(x) - u(y)|^{p- \eps}}{|x-y|^{p -  \eps }}  \frac{\eps}{|x-y|^{N - \frac{\eps \gamma}{p}}}d x \, dy \\[6pt]
\ge \frac{1}{p-\eps} \frac{1}{|2m|^{2\eps}} \mathop{\mathop{\iint}_{x \in B_R; |x-y| < r}}_{\frac{|u(x) - u(y)|}{|x - y|^{1+ \frac{\gamma}{p}}} > 1}       \frac{|u(x) - u(y)|^{p+\eps}}{|x-y|^{p+\eps}}  \frac{\eps}{|x-y|^{N - \frac{\eps \gamma}{p} - 2 \eps}}d x \, dy. 
\end{multline}
Combining \eqref{thm-neg-p1}, and \eqref{thm-neg-p2} yields 
\begin{multline} \label{thm-neg-p3}
\int_1^\infty \frac{\eps}{\lambda^{1 + \eps}} \Phi_{\lambda} (u) \, d \lambda 
\ge \frac{1}{p-\eps} \frac{1}{|2m|^{2\eps}} \mathop{\mathop{\iint}_{x \in B_R; |x-y| < r}}_{\frac{|u(x) - u(y)|}{|x - y|^{1+ \frac{\gamma}{p}}} > 1}      \frac{|u(x) - u(y)|^{p+\eps}}{|x-y|^{p+\eps}}  \frac{\eps}{|x-y|^{N - \frac{\eps \gamma}{p} - 2 \eps}}d x \, dy \\[6pt]
- \frac{\eps}{p-\eps}   \mathop{\mathop{\iint}_{x \in B_R; |x-y| < r}}_{\frac{|u(x) - u(y)|}{|x - y|^{1+ \frac{\gamma}{p}}} > 1}  |x-y|^{-N + \gamma} \, dx \, dy.
\end{multline}

Similarly, by \eqref{thm-neg-p0-1}, it holds 
\be \label{thm-neg-p4}
\limsup_{\lambda \to 0_+} \Phi_\lambda (u) \ge \limsup_{\eps \to 0_+} \int_0^1 \frac{\eps}{\lambda^{1 - \eps}} \Phi_{\lambda} (u) \, d \lambda. 
\ee
We have 
\begin{multline} \label{thm-neg-p5}
\int_0^1 \frac{\eps}{\lambda^{1 - \eps}} \Phi_{\lambda} (u) \, d \lambda \ge \frac{1}{p + \eps} \mathop{\mathop{\iint}_{x \in B_R; |x-y| < r}}_{\frac{|u(x) - u(y)|}{|x - y|^{1+ \frac{\gamma}{p}}} \le 1}     \frac{|u(x) - u(y)|^{p+\eps}}{|x-y|^{p + \eps}} \frac{\eps}{|x-y|^{N + \frac{\eps \gamma}{p}}}  d x \, dy \\[6pt]
 +  \frac{\eps}{p+ \eps}   \mathop{\mathop{\iint}_{x \in B_R; |x-y| < r}}_{\frac{|u(x) - u(y)|}{|x - y|^{1+ \frac{\gamma}{p}}} > 1}  |x-y|^{-N + \gamma} \, dx \, dy. 
\end{multline}

Combining \eqref{thm-neg-p3} and \eqref{thm-neg-p5} yields, since $\frac{\eps \gamma }{p} > - \frac{\eps \gamma }{p} - 2 \eps$ and $ r< 1$,  
\begin{multline}\label{thm-neg-p6}
\int_1^\infty \frac{\eps}{\lambda^{1 + \eps}} \Phi_{\lambda} (u) \, d \lambda  + \frac{p+\eps}{p-\eps} \int_0^1 \frac{\eps}{\lambda^{1 - \eps}} \Phi_{\lambda} (u) \, d \lambda  \\[6pt] 
\ge \frac{1}{p-\eps} \frac{1}{|2m|^{2\eps}} \mathop{\iint}_{x \in B_R; |x-y| < r}    \frac{|u(x) - u(y)|^{p+\eps}}{|x-y|^{p + \eps}} \frac{\eps}{|x-y|^{N - \frac{\eps \gamma}{p} - 2 \eps}}  d x \, dy.
\end{multline}

Let $(\varphi_k)_{k \ge 1}$ be a smooth non-negative sequence of approximations to the identity such that $\supp \varphi_k \subset B_{1/k}$. Set 
$$
u_k = \varphi_k  * u \mbox{ for } k \ge 1. 
$$
We have, by Jensen's inequality,  
\begin{multline}\label{thm-neg-p7}
 \frac{1}{p-\eps} \frac{1}{|2m|^{2 \eps}}\mathop{\iint}_{x \in B_R; |x-y| < r}    \frac{|u(x) - u(y)|^{p + \eps}}{|x-y|^{p + \eps}} \frac{\eps}{|x-y|^{N - \frac{\eps \gamma}{p} - 2 \eps}} d x \, dy \\[6pt] \ge   \frac{1}{p-\eps} \frac{1}{|2m|^{2 \eps}} \mathop{\iint}_{x \in B_{R-1}; |x-y| < r}    \frac{|u_k(x) - u_k(y)|^{p + \eps}}{|x-y|^{p + \eps}} \frac{\eps}{|x-y|^{N - \frac{\eps \gamma}{p} - 2 \eps}} d x \, dy.
\end{multline}
Since 
\be
\liminf_{\eps \to 0_+} \mathop{\iint}_{x \in B_{R-1}; |x-y| < r}    \frac{|u_k(x) - u_k(y)|^{p + \eps}}{|x-y|^{p + \eps}} \frac{\eps}{|x-y|^{N - \frac{\eps \gamma}{p} - 2 \eps}} d x \, dy \ge \frac{K_{N,p}}{\gamma  + 2p}\int_{B_{R-1}} |\nabla u_k|^p \, dx, 
\ee
it follows from \eqref{thm-neg-p6} and \eqref{thm-neg-p7} that 
\be
\liminf_{\eps \to 0_+} \left( \int_1^\infty \frac{\eps}{\lambda^{1 + \eps}} \Phi_{\lambda} (u) \, d \lambda  + \int_0^1 \frac{\eps}{\lambda^{1 - \eps}} \Phi_{\lambda} (u) \, d \lambda \right) \ge 
\frac{K_{N,p}}{\gamma  + 2p}\int_{B_{R-1}} |\nabla u_k|^p \, dx. 
\ee
Since $R>2$ is arbitrary, we reach 
\be
\limsup_{\lambda \to 0_+} \Phi_\lambda(u) + \limsup_{\lambda \to + \infty} \Phi_\lambda(u)  \ge \frac{K_{N,p}}{\gamma + 2p}  \int_{\mR^N} |\nabla u_k|^p \, dx. 
\ee
Letting $k \to + \infty$, we reach the conclusion when $u \in L^\infty(\mR^N)$. 

The proof in the general case follows from the case where $u \in L^\infty(\mR^N)$ as in the proof of \Cref{thm-I} after noting that 
$$
\Phi_\lambda (u) \ge \Phi_\lambda (u_A), 
$$
where 
$$
u_A = \min\{ \max\{u, - A\}, A\}. 
$$
The details are omitted.
\qed
}

\section{Proof of \Cref{pro-CE}} \label{sect-pro-CE}

This section consisting of two subsections is devoted to the proof of \Cref{pro-CE}. In the first subsection, we present state and prove a lemma which is used in the proof of \Cref{pro-CE}. The proof of \Cref{pro-CE} is given in the second subsection.

\subsection{A useful lemma}

Here is the result of this section. 

\begin{lemma}\label{lem1} Let $1 \le p < + \infty$, $\beta \in \mR$,  and let $g$ be a measurable function on $\mR^N$. There exist positive constants $c_1, c_2$ independent of $g$ such that, for $\gamma \in \mR \setminus \{0 \}$, 
\begin{multline} \label{lem1-cl1}
\mathop{\int_{\mR^{N-1}} \int_{\mR} \int_{\mR}}_{\frac{|g(x_i', x_i) -
g(x_i', y_i)|}{|x_i - y_i|^{\beta}} > c_1 \lambda} |x_i - y_i|^{-1 + \gamma} \, dx_i \,
dy_i \, dx_i' \\[6pt]
\le c_2  \mathop{\int_{\mR^N} \int_{\mR^N}}_{\frac{|g(x)
- g(y)|}{|x-y|^\beta} > \lambda} |x-y|^{-N+\gamma} \, dx \, dy \quad
\forall \, \lambda>0, \mbox{ for } 1 \le i \le N,  
\end{multline}
and, for $\gamma < 1$ \hmr{ and $\beta \ge 0$}, 
\begin{equation}\label{lem1-cl2}
\mathop{\int_{\mR^N} \int_{\mR^N}}_{\frac{|g(x)
- g(y)|}{|x-y|^\beta} > c_1 \lambda} |x-y|^{-N+\gamma} \, dx \, dy \le c_2  \sum_{i=1}^N \mathop{\int_{\mR^{N-1}} \int_{\mR} \int_{\mR}}_{\frac{|g(x_i', x_i) - g(x_i', y_i)|}{|x_i - y_i|^{\beta}} > \lambda} |x_i - y_i|^{-1 + \gamma} \, dx_i \,
dy_i \, dx_i', \forall \, \lambda>0. 
\end{equation}
\end{lemma}

Here and in what follows, for $1 \le i \le N$ and for a mesurable function $f$ defined in $\mR^N$, we denote
$$
f(x_i', x_i) = f (x) \mbox{ for } x = (x_1, \cdots, x_N) \in \mR^N
$$
where 
$$
x_i' = (x_1, \cdots, x_{i-1}, x_{i+1}, \cdots, x_N) \in \mR^{N-1}. 
$$

\begin{proof} We begin with the proof of \eqref{lem1-cl1}. For notational ease, we only consider the case $i=N$ and denote $x_N'= x'$ in this case. We have 
\begin{multline} \label{lem1-cl1-p1}
\mathop{\int_{\mR^{N-1}} \int_{\mR} \int_{\mR}}_{\frac{|g(x', x_N) -
g(x', y_N)|}{|x_N- y_N|^{\beta}} > c_1 \lambda} |x_N - y_N|^{-1 + \gamma} \, dx_N  \,
dy_N \, dx' 
\\[6pt]
=  \mathop{\int_{\mR^{N-1}} \int_{\mR} \int_{\mR}}_{\frac{|g(x', x_N) -
g(x', y_N)|}{|x_N- y_N|^{\beta}} > c_1 \lambda} \fint_{B_{
\mR^{N-1}}(x', |y_N - x_N |/2)} |x_N - y_N|^{-1 + \gamma} \, dz' \, dx_N  \,
dy_N \, dx'. 
\end{multline}
Here, for $z' \in \mR^{N-1}$ and $r>0$, we denote $B_{\mR^{N-1}} (z', r)$ the open ball in $\mR^{N-1}$ centered at $z' \in \mR^{N-1}$ and of radius $r>0$, and $\fint_{B_{\mR^{N-1}}(z', r)} = \frac{1}{|B_{\mR^{N-1}(z', r)}|} \int_{B_{\mR^{N-1}}(z', r)}$.

Since $|g(x', x_N) - g(x', y_N)| \le |g(x', x_N) - g(z, \frac{x_{N} + y_N}{2})| + |g(z, \frac{x_{N} + y_N}{2}) -
g(x', y_N)|$, we obtain
\begin{multline}\label{lem1-cl1-p2}
\mathop{\int_{\mR^{N-1}} \int_{\mR} \int_{\mR}}_{\frac{|g(x', x_N) -
g(x', y_N)|}{|x_N- y_N|^{\beta}} > c_1 \lambda} \fint_{B_{
\mR^{N-1}}(x', |y_N - x_N |/2)} |x_N - y_N|^{-1 + \gamma} \, dz' \, dx_N  \,
dy_N \, dx'
\\[6pt]
\le \mathop{\int_{\mR^{N-1}} \int_{\mR} \int_{\mR}}_{\frac{|g(x', y_N) - g(z', \frac{x_{N} + y_N}{2})|}{|x_N- y_N|^{\beta}} > c_1\lambda/2} \fint_{B_{
\mR^{N-1}(x', |y_N - x_N |/2)}} |x_N- y_N|^{-1 + \gamma} \, dz' \, dx_N  \,
dy_N \, dx' \\[6pt]
\hmr{+} \mathop{\int_{\mR^{N-1}} \int_{\mR} \int_{\mR}}_{\frac{|g(x', x_N) - g(z', \frac{x_{N} + y_N}{2})|}{|x_N-  y_N|^{\beta}} > c_1\lambda/2} \fint_{B_{
\mR^{N-1}(x', |y_N - x_N |/2)}} |x_N- y_N|^{-1 + \gamma} \, dz' \, dx_N  \,
dy_N \, dx'. 
\end{multline}
We have, by a change of variables, for $c_1$ sufficiently large, 
\begin{multline}\label{lem1-cl1-p3}
\mathop{\int_{\mR^{N-1}} \int_{\mR} \int_{\mR}}_{\frac{|g(x', x_N) - g(z', \frac{x_{N} + y_N}{2})|}{|x_N -  y_N|^{\beta}} > c_1\lambda/2} \fint_{B_{
\mR^{N-1}}(x', |y_N - x_N |/2)} |x_N- y_N|^{-1 + \gamma} \, dz' \, dx_N  \,
dy_N \, dx' \\[6pt]
\le C \mathop{\iint_{\mR^{N} \times \mR^N}}_{\frac{|g(x) - g(y)|}{|x-y|^{\beta}} > \lambda}  |x-y|^{-N+ \gamma}\, dx \, dy,
\end{multline}
and 
\begin{multline}\label{lem1-cl1-p4}
\mathop{\int_{\mR^{N-1}} \int_{\mR} \int_{\mR}}_{\frac{|g(x', y_N) - g(z', \frac{x_{N} + y_N}{2})|}{|x_N- y_N|^{\beta}} > c_1 \lambda/2} \fint_{B_{
\mR^{N-1}}(x', |y_N - x_N |/2)} |x_N- y_N|^{-1 + \gamma} \, dz' \, dx_N  \,
dy_N \, dx'  \\[6pt] \le C \mathop{\iint_{\mR^{N} \times \mR^N}}_{\frac{|g(x) - g(y)|}{|x-y|^{\beta}} > \lambda}  |x-y|^{-N+ \gamma}\, dx \, dy.
\end{multline}
Combining \eqref{lem1-cl1-p1}, \eqref{lem1-cl1-p2}, \eqref{lem1-cl1-p3}, and \eqref{lem1-cl1-p4} yields \eqref{lem1-cl1}.    

We next give the proof of \eqref{lem1-cl2}. For notational ease, we only consider the case $N=2$. The general case follows similarly. Since 
$|g(x_1, x_2) - g(y_1, y_2)| \le |g(x_1, x_2) - g(x_1, y_2)| + |g(x_1, y_2) - g(x_2, y_2)|$, \hmr{ and $\beta \ge 0$,} we derive that, for $c_1$ sufficiently large, 
\begin{multline} \label{lem1-cl2-p1}
\mathop{\iint_{\mR^{2} \times \mR^2}}_{\frac{|g(x_1, x_2) - g(y_1, y_2)|}{|x-y|^{\beta}} > c_1 \lambda}  |x-y|^{-2+ \gamma}\, dx \, dy \\[6pt]
\le \mathop{\iint_{\mR^{2} \times \mR^2}}_{\frac{|g(x_1, x_2) - g(x_1, y_2)|}{|x_2-y_2|^{\beta}} > \lambda}  |x-y|^{-2+ \gamma}\, dx \, dy + \mathop{\iint_{\mR^{2} \times \mR^2}}_{\frac{|g(x_1, y_2) - g(y_1, y_2)|}{|x_1 - y_1|^{\beta}} > \lambda}  |x-y|^{-2+ \gamma}\, dx \, dy.
\end{multline}
Since $\gamma < 1$, it follows that 
\be \label{lem1-cl2-p2} 
\mathop{\iint_{\mR^{2} \times \mR^2}}_{\frac{|g(x_1, x_2) - g(x_1, y_2)|}{|x_2-y_2|^{\beta}} > \lambda}  |x-y|^{-2+ \gamma}\, dx \, dy \\[6pt]
 \le C \int_{\mR} \mathop{\iint_{\mR \times \mR}}_{\frac{|g(x_1, x_2) - g(x_1, y_2)|}{|x_2-y_2|^{\beta}} > \lambda}  |x_2-y_2|^{-1+ \gamma}\, dx_1 \, dx_2 \, dy_2, 
\ee
and
\be \label{lem1-cl2-p3}
\mathop{\iint_{\mR^{2} \times \mR^2}}_{\frac{|g(x_1, y_2) - g(y_1, y_2)|}{|x_1 - y_1|^{\beta}} > \lambda}  |x-y|^{-2+ \gamma}\, dx \, dy \le C \int_{\mR} \mathop{\iint_{\mR^{2} \times \mR^2}}_{\frac{|g(x_1, y_2) - g(y_1, y_2)|}{|x_1 - y_1|^{\beta}} > \lambda}  |x_1-y_1|^{-1+ \gamma}\, dx_1 \, dy_1 \, d y_2. 
\ee 
Combining \eqref{lem1-cl2-p1}, \eqref{lem1-cl2-p2}, and \eqref{lem1-cl2-p3} yields \eqref{lem1-cl2}.   

\medskip 
The proof is complete. 
\end{proof}

\begin{remark} \rm The method used to
prove Lemma~\ref{lem1} \hmr{appeared} in the theory of fractional Sobolev spaces (see, e.g., \cite[Chapter 7]{Adams-75}). \hmr{In this context, it is due to Besov.} 
\end{remark}

\subsection{Proof of \Cref{pro-CE}}
Set
$$
\beta: = 1 + \gamma/p. 
$$
Then $\beta > 0$ since $-1 < \gamma < 0$.  Consider
$$
u = \mathds{1}_{Q} \mbox{ where } Q = (0, 1)^N \subset \mR^N. 
$$
By \Cref{lem1}, there exist two positive constants $c_1, c_2$ such that, for all $\lambda > 0$,  
\begin{multline}\label{pro-CE-p1}
\mathop{\int_{\mR^N} \int_{\mR^N}}_{\frac{|u(x)
- u(y)|}{|x-y|^\beta} > c_1 \lambda} |x-y|^{-N+\gamma} \, dx \, dy \le c_2  \sum_{i=1}^N \mathop{\int_{\mR^{N-1}} \int_{\mR} \int_{\mR}}_{\frac{|u(x_i', x_i) - u(x_i', y_i)|}{|x_i - y_i|^{\beta}} >  \lambda} |x_i - y_i|^{-1 + \gamma} \, dx_i \,
dy_i \, dx_i' \\[6pt]
=  c_2 N \mathop{\int_{\mR^{N-1}} \int_{\mR} \int_{\mR}}_{\frac{|u(x', x_N) - u(x', y_N)|}{|x_N - y_N|^{\beta}} >  \lambda} |x_N - y_N|^{-1 + \gamma} \, dx_N \,
dy_N \, dx', 
\end{multline}
where we denote $x_N'$ by $x'$ for notational ease. Since
\begin{multline*}
\left\{ (x', x_N, y_N) \in \mR^{N-1} \times \mR \times \mR ; \frac{|u(x', x_N) - u(x', y_N)|}{|x_N - y_N|^{\beta}} > \lambda \right\} \\[6pt]
=\left\{ (x', x_N, y_N) \in (0, 1)^{N-1} \times \hmr{\Big( \big( (0, 1) \times (\mR \setminus (0, 1)) \big)\cup \big( (\mR \setminus (0, 1)) \times (0, 1) \big) \Big)} ; |x_N - y_N|^{\beta} < 1/ \lambda \right\},
\end{multline*}
it follows that, for $\lambda^{1/\beta} < 1/4$,  
\be \label{pro-CE-p2}
\lambda^p \mathop{\int_{\mR^{N-1}} \int_{\mR} \int_{\mR}}_{\frac{|u(x', x_N) - u(x', y_N)|}{|x_N - y_N|^{\beta}} >  \lambda} |x_N - y_N|^{-1 + \gamma} \, dx_N \,
dy_N \, dx' \le C \lambda^p \int_0^1 x_N^{\gamma} \, dx_N \le C_{\gamma} \lambda^p.  
\ee
The conclusion now follows from \eqref{pro-CE-p1} and \eqref{pro-CE-p2}. \qed

\section{Proof of \Cref{thm-neg-1p}} \label{sect-thm-neg-1p}

\subsection{A fundamental lemma}

We follow the approach of Bourgain and Nguyen \cite{BN-06}.  The following lemma is the key ingredient.

\begin{lemma}\label{lem-fund} Let $1\le p < + \infty$, $-p \le \gamma \le -1$, and let $f$ be a measurable function on a
bounded open nonempty interval $I$. Then
\begin{equation}\label{lem-fund-cl}
\dsp \mathop{\liminf}_{\lambda \to 0_+}
\Phi_{\hmr{\lambda}, I} (f)
 \ge c \frac{1}{| I|^{p-1}}(\esssup f - \essinf f)^{p},
\end{equation}
where $c=c_p$ is a positive constant depending only on $p$.
\end{lemma}

\begin{proof} Set 
$$
\beta = 1 + \gamma / p. 
$$
Then 
\be \label{lem-fund-beta}
0 \le \beta < 1. 
\ee
In what follows, we assume that 
\be \label{lem-fund-p0}
\mathop{\liminf}_{\lambda \to 0_+}
\Phi_{\hmr{\lambda}, I} (u) < + \infty 
\ee
since there is nothing to prove otherwise. 

\medskip 
The proof is now divided into two steps. 

\medskip 
\noindent \underline{Step 1.} Proof of \eqref{lem-fund-cl}  for $f \in
L^\infty(I)$.

\medskip 
By rescaling, we may assume $ I= (0,1)$. Denote $s_+ = \esssup f$ and  $s_- = \essinf f$. Rescaling $f$, one
may also assume
\begin{equation}\label{varf}
s_+ - s_- =1
\end{equation}
(unless $f$ is constant on $ I$ and there is nothing to prove in this case).

Take $0< \delta < 1$ small enough to ensure that there are
(density) points $t_+$, $t_- \in [40 \delta, 1- 40\delta] \subset
[0,1]$ with
\begin{equation}\label{pointt}
\left\{\begin{array}{l} \dsp \left| \left[ t_+ - \tau, t_+ +
\tau\right] \cap \left[f
> \frac{3}{4}s_+ + \frac{1}{4}s_- \right]\right| > \frac{9}{5}
\tau,\\[6pt]
\dsp \left| \left[ t_- - \tau, t_- + \tau\right] \cap \left[f <
\frac{3}{4}s_- + \frac{1}{4}s_+ \right]\right| > \frac{9}{5} \tau,
\end{array}\right. \quad \forall  \, 0< \tau < 40 \delta.
\end{equation}
Take $K \in \ZZ_+$ such that $\delta < 2^{-K} \le 5 \delta/4$ and
denote
\begin{equation*}
J =\left\{ j \in \ZZ_+; \frac{3}{4}s_- + \frac{1}{4}s_+ < j 2^{-K}
< \frac{3}{4} s_+ + \frac{1}{4} s_- \right\}.
\end{equation*}
Then
\begin{equation}\label{cardJ}
|J| \ge 2^{K-1} - 2  \approx \frac{1}{\delta}.
\end{equation}
For each $j$, define the following sets
\begin{equation*}
\begin{array}{c}
A_j = \left\{ x \in [0,1]; \, (j-1)2^{-K} \le f(x) < j
2^{-K}\right\}, \quad  B_j = \bigcup_{j'< j}A_{j'} \quad
\mbox{and} \quad C_j = \bigcup_{j'>j}A_{j'},
\end{array}
\end{equation*}
so that $B_j \times C_j \subset \left[|f(x)-f(y)| \ge
2^{-K}\right] \subset \left[|f(x)-f(y)| > \delta\right]$.

Since the sets $A_j$ are disjoint, it follows from \eqref{cardJ}
that
\begin{equation}\label{cardG}
\mathrm{card}(G) \ge 2^{K-2} - 3 \approx \frac{1}{\delta},
\end{equation}
where $G$ is defined by
\begin{equation*}
G=\{j \in J; \, |A_j| < 2^{-K+2}\}.
\end{equation*}

For each $j \in G$, set $\lambda_{1,j} = |A_j|$ and  consider the
function $\psi_1(t)$ defined as follows
\begin{equation*}
\psi_1(t)=\left|[t- 4 \lambda_{1,j}, t + 4 \lambda_{1,j}] \cap B_j
\right|, \quad \forall \, t \in [40 \delta , 1- 40 \delta].
\end{equation*}
Then, from \eqref{pointt}, $\psi_1(t_+) < 4 \lambda_{1,j}$ and
$\psi_1(t_-) > 4 \lambda_{1,j}$. Hence, since $\psi_1$ is a
continuous function on the interval $[40 \delta, 1- 40 \delta ]$
containing the two points $t_+$ and $t_-$, there exists $t_{1,j}
\in [40 \delta, 1 - 40 \delta ]$ such that
\begin{equation}\label{deft1}
\psi_1(t_{1,j})= 4 \lambda_{1,j}.
\end{equation}
We have 
$$  
\mathop{\mathop{\iint}_{ I \times I}}_{|f(x)-
f(y)|> \delta} |x-y|^{-2} \, dx \, dy< + \infty
$$ 
by \eqref{lem-fund-p0} and the fact that $\beta \ge 0$ and $\gamma \le - 1$.  It
follows that $\left|[t_{1,j}- 4 \lambda_{1,j}, t_{1,j} + 4
\lambda_{1,j}] \cap A_j \right|>0$. Indeed, suppose $\left|[t_{1,j}- 4 \lambda_{1,j}, t_{1,j} + 4
\lambda_{1,j}] \cap A_j \right|=0$. Then
\begin{equation*}
\mathop{\mathop{\iint}_{x \in [t_{1,j}- 4 \lambda_{1,j}, t_{1,j} +
4 \lambda_{1,j}] \cap B_j}}_{y \in [t_{1,j}- 4 \lambda_{1,j},
t_{1,j} + 4 \lambda_{1,j}] \setminus B_j}\frac{1}{|x-y|^{2}} \, dx
\, dy  \le \mathop{\mathop{\iint}_{ I \times I}}_{|f(x)- f(y)|>
\delta} \frac{1}{|x-y|^{2}} \, dx \, dy < + \infty.
\end{equation*}
Hence $|[t_{1,j}- 4 \lambda_{1,j}, t_{1,j} + 4 \lambda_{1,j}] \cap
B_j| =0$ or $|[t_{1,j}- 4 \lambda_{1,j}, t_{1,j} + 4
\lambda_{1,j}] \setminus B_j| =0$. This is a
contradiction since $\psi_1(t_{1,j})=|[t_{1,j}- 4 \lambda_{1,j},
t_{1,j} + 4 \lambda_{1,j}] \cap B_j|=  4
\lambda_{1,j}$ (see \eqref{deft1}).

If $ \dsp \left|[t_{1,j}- 4 \lambda_{1,j}, t_{1,j} + 4
\lambda_{1,j}] \cap A_j \right| < \lambda_{1,j}/4,$ then take
$\lambda_{2,j}>0$ such that $\lambda_{1,j}/\lambda_{2,j} \in
\ZZ_+$ and
\begin{equation*}
\frac{\lambda_{2,j}}{2} < \left|[t_{1,j}- 4 \lambda_{1,j}, t_{1,j}
+ 4 \lambda_{1,j}] \cap A_j \right| \le \lambda_{2,j}.
\end{equation*}
Since $\dsp \left|[t_{1,j}- 4 \lambda_{1,j}, t_{1,j} + 4
\lambda_{1,j}] \cap A_j \right| < \lambda_{1,j}/4$, we infer that
$ \dsp
\lambda_{2,j} \le \lambda_{1,j}/2$.  

Set $E_{2,j} = [t_{1,j} - 4 \lambda_{1,j} + 4 \lambda_{2,j},
t_{1,j} + 4 \lambda_{1,j} - 4 \lambda_{2,j}]$ and consider the
function $\psi_2(t)$ defined as follows
\begin{equation*}
\psi_2(t)=\left|[t- 4 \lambda_{2,j}, t + 4 \lambda_{2,j}] \cap B_j
\right|, \quad \forall \, t \in E_{2,j}.
\end{equation*}

We claim that there exists $t_{2,j} \in E_{2,j}$ such that
$\psi_2(t_{2,j}) = 4 \lambda_{2,j}$.

To see this, we argue by contradiction. Suppose that $\psi_2(t)
\neq 4 \lambda_{2,j}$, for all $t \in E_{2,j}$. Since $\psi_2$ is
a continuous function on $E_{2,j}$, we assume as well that
$\psi_2(t)< 4 \lambda_{2,j}$, for all $t \in E_{2,j}$. Since
$\lambda_{1,j}/\lambda_{2,j} \in \ZZ_+$, it follows that
$\psi_1(t_{1,j}) < 4
\lambda_{1,j}$, hence we have a contradiction to \eqref{deft1}. Thus the claim is proved. 

It is clear that
\begin{equation*}
\mathop{\mathop{\iint}_{[t_{2,j}- 4 \lambda_{2,j}, t_{2,j} + 4
\lambda_{2,j}]^2}}_{|f(x)- f(y)|> \delta} \frac{1}{|x-y|^{2}} \,
dx \, dy \ge \mathop{\mathop{\iint}_{[t_{2,j}- 4 \lambda_{2,j},
t_{2,j} + 4 \lambda_{2,j}]^2} }_{x \in B_j; \, y \in C_j}
\frac{1}{|x-y|^{2}} \, dx \, dy \sge 1.
\end{equation*}

If $ \dsp \left|[t_{2,j}- 4 \lambda_{2,j}, t_{2,j} + 4
\lambda_{2,j}] \cap A_j \right| < \lambda_{2,j}/4, $
then take $\lambda_{3,j}$ ($\lambda_{3,j} \le \dsp \lambda_{2,j}/2$) and  $t_{3,j}$, etc. 
On the other hand, since $\dsp \mathop{\mathop{\iint}_{ I \times
 I}}_{|f(x)- f(y)|> \delta} \frac{1}{|x-y|^{2}} \, dx \, dy< +
\infty, $ we have
\begin{equation}\label{proborne}
\mathop{\limsup_{\mathrm{diam}(Q) \lr 0}}_{Q \mathrm{ :\, an \;
interval\; of \;}  I} \mathop{\mathop{\iint}_{Q \times Q}}_{|f(x)-
f(y)|> \delta} \frac{1}{|x-y|^{2}} \, dx \, dy =0.
\end{equation}
Thus, from \eqref{proborne} and the construction of $t_{k,j}$ and
$\lambda_{k,j}$, there exist $t_j \in [40 \delta, 1- 40 \delta]$
and $\lambda_j > 0$ ($t_j = t_{k, j}$ and $\lambda_j =
\lambda_{k,j}$ for some $k$) such that
\begin{equation}\label{measAj}
a) \, \left|[t_j- 4 \lambda_j, t_j + 4 \lambda_j] \cap B_j \right|
= 4 \lambda_j \quad \mbox{and} \quad b) \, \frac{\lambda_j}{4} \le
\left|[t_j- 4 \lambda_j, t_j + 4 \lambda_j] \cap A_j \right| \le
\lambda_j.
\end{equation}

Set $\lambda =    \inf_{j \in G} \lambda_j$ ($\lambda>0$ since $G$
is finite). Suppose $G=\bigcup_{i=1}^{n} I_m$, where $I_m$ is
defined as follows
\begin{equation*}
I_{m} = \Big\{ j \in G; \, 2^{m-1} \lambda \le \lambda_j < 2^{m}
\lambda \Big\}, \quad \forall \, m \ge 1.
\end{equation*}
Then it follows from \eqref{cardG} that
\begin{equation}\label{cardJ1m}
\sum_{m=1}^{ n }\card{I_{m}} \gtrsim \frac{1}{\delta}.
\end{equation}
For each $m$ ($1 \le m \le n$), since $A_j \cap A_k = \O$ for $j
\neq k$, it follows from (\ref{measAj}-b) that there exists $J_m
\subset I_m$ such that
\begin{equation}\label{defJm0}
 a) \, \card{J_m} \gtrsim \card{I_m} \quad \mbox{and} \quad
b)\,  |t_i - t_j| > 2^{m + 3} \lambda, \quad  \, \forall \, i, j
\in J_m.
\end{equation}
Then, from (\ref{defJm0}-b) and the definition of $I_m$,
\begin{equation}\label{defJm}
[t_i - 4 \lambda_i, t_i + 4 \lambda_i] \cap [t_j - 4 \lambda_j,
t_j + 4 \lambda_j] = \O, \quad \forall \, i, j \in J_m.
\end{equation}
Set $U_0 : = \O$ and
\begin{equation*}
\left\{
\begin{array}{l}  L_m = \dsp \big\{ j \in J_m; \, \big|[t_j - 4
\lambda_j, t_j + 4 \lambda_j] \setminus U_{m-1} \big| \ge 6
\lambda_j \big\},\\[6pt]
 U_m =  \big(\bigcup_{j \in L_m}[t_j - 4 \lambda_j, t_j + 4
\lambda_j]\big) \cup U_{m-1},\\[6pt]
a_m = \card{J_m} \quad \mbox{and}\quad b_m = \card{L_m},
\end{array} \right. \quad \mbox{for } m =1, 2, \dots, n.
\end{equation*}
From \eqref{defJm} and the definitions of $J_m$ and $L_m$,
\begin{equation*}
\frac{1}{4}2^{m -1} (a_m - b_m) \le \sum_{i=1}^{m-1} 2^{i} b_i
\end{equation*}
which shows that
\begin{equation*}
a_m \le b_m + 8 \sum_{i=1}^{m-1} 2^{(i-m)} b_i.
\end{equation*}
Consequently,
\begin{align*}
\sum_{m=1}^n a_m & \le \sum_{m=1}^n b_m + 8  \sum_{m=1}^n
\sum_{i=1}^{m-1} 2^{(i-m)} b_i =  \sum_{m=1}^n b_m + 8
\sum_{i=1}^n b_i \sum_{m=i+1}^{n}2^{(i-m)}.
\end{align*}
Since $\sum_{i=1}^\infty 2^{-i} =1$, it follows from from
\eqref{cardJ1m} and (\ref{defJm0}-a) that
\begin{equation*}
\sum_{m=1}^n b_m \ge \frac{1}{9}\sum_{m=1}^n a_m \sge
\frac{1}{\delta}.
\end{equation*}

Take $\lambda$ such that 
$$
\delta \ge \lambda (8 \delta)^\beta \ge \delta/2. 
$$
 We then have, since $\gamma \le -1$,  
\begin{multline*}
\mathop{\mathop{\iint}_{ I \times  I}}_{\frac{|f(x)- f(y)|}{|x-y|^\beta}> \lambda}
\lambda^p |x-y|^{-1 + \gamma}\, dx \, dy \ge \sum_{m=1}^n \sum_{j \in L_m}
\mathop{\mathop{\iint}_{\left([t_j - 4 \lambda_j, t_j + 4
\lambda_j] \setminus U_{m-1}\right)^2}}_{x \in B_j, \, y \in C_j}
\lambda^p |x-y|^{-1 + \gamma}\, dx \, dy \\[6pt]
\sge \sum_{m=1}^n
b_m \lambda^p \delta^{1+\gamma}  \sge \lambda^p \delta^{\gamma} \sge 1. 
\end{multline*}
which yields \eqref{lem-fund-cl} by \eqref{lem-fund-beta}. 

\underline{Step 2.} Proof of \Cref{lem-fund-cl} in the general case.

For $A> 0$, denote 
\begin{equation*}
f_A = (f \vee (-A))\wedge A,
\end{equation*}
then
\begin{equation*}
|f_A(x) - f_A(y)| \le |f(x) - f(y)| \mbox{ for } (x, y) \in I^2. 
\end{equation*}
Applying the result of Step 1 to the sequence $f_A$ and letting $A$
goes to infinity, we deduce that \eqref{lem-fund-cl} holds for any
measurable function $f$ on $ I$. 

\medskip 
The proof is complete. 
\end{proof}

\subsection{Proof of \Cref{thm-neg-1p}} 

\noindent {\it Step 1:} Proof of \Cref{thm-neg-1p} when $N=1$. 

\medskip Set
$\tau_{h}(u)(x) = \dsp \frac{u(x+h) - u (x)}{h}, \, \forall \, x
\in
\mR$, $0< h< 1$. For each $m \ge 2$, take $K\in \mR_+$ such that $Kh \ge m$, then
\begin{equation*}
\int_{-m}^m |\tau_h(u)(x)|^p \, dx \le \sum_{k = -K}^{K}\int_{k
h}^{(k+1)h } |\tau_h(u)(x)|^p \, dx.
\end{equation*}
Thus, since
\begin{equation*}
\int_a^{a + h } |\tau_{h}(u)(x)|^p \, dx \le \int_a^{a + h
}\frac{1}{h^p} |\mathop{\mathrm{ess \; sup} }_{x \in (a, a + 2h )}
u - \mathop{\mathrm{ess \; inf} }_{x \in (a, a + 2h ) } u|^p \,
dx,
\end{equation*}
it follows from \Cref{lem-fund} that, for some constant
$c=c_p>0$,
\begin{equation}\label{compactbound}
\int_{-m}^m |\tau_h(u)(x)|^p \, dx \le c  \liminf_{\lambda \to 0} \Phi_\lambda (u). 
\end{equation}
Since $m \ge 2$ is arbitrary, \eqref{compactbound} shows that
\begin{equation}\label{estth}
\int_{\mR} |\tau_h(u)(x)|^p \, dx \le c  \liminf_{\lambda \to 0} \Phi_\lambda (u). 
\end{equation}
Since \eqref{estth} holds for all $0< h < 1$, it follows that $\nabla 
\in L^p(\mR)$ for $p>1$ (see e.g. \cite[Chapter 8]{Brezis-FA}) and $\nabla u \in \cM$ (see, e.g., \cite{EGMeasure}). 
Moreover, \eqref{thm-neg-1p-cl} holds. 

\medskip 
\noindent {\it Step 2:} Proof of \Cref{thm-neg-1p} when $N \ge 2$. The result in this case is a consequence of \Cref{lem-fund}, \Cref{lem1}, and \cite[Proposition \hmr{3}]{Ng-11}. As in the proof of the case $N=1$, one can show that $\mbox{ess} V(u, j)$ is bounded by $C \liminf_{\lambda \to 0} \Phi_\lambda (u)$ for $1 \le j \le N$ where $\mbox{ess} V (u, j)$ is the essential variation of $u$ in the $j$-th direction whose definition is given in \cite[Definition 5]{Ng-11}. The details are omitted. \qed

\section{Gamma-convergence} \label{sect-Gamma}

In this section, we discuss the $\Gamma$-convergence of the families of nonlocal functionals mentioned in the introduction.  Concerning the BBM formula, the following result was established by Ponce \cite{Ponce-04}. 

\begin{theorem}[Ponce] \label{thm-BBM-Gamma} Let $N \ge 1$, $1 \le  p < + \infty$, and let $(\rho_n)_{n \ge 1}$ be a sequence of non-negative  mollifiers.  Denote 
$$
J_n(u): = \mathop{\iint}_{\mR^N \times \mR^N} \frac{|u(x)
- u(y)|^p}{|x-y|^p} \rho_n(|x-y|)\, dx \, dy. 
$$
Then 
$$
J_n \mbox{ $\Gamma$-converges to } K_{N, p} \Phi \mbox{ in } L^p(\mR^N) \mbox{ as } n \to + \infty. 
$$
\end{theorem}

To establish \Cref{thm-BBM-Gamma}, one needs to show that, given $u \in L^p(\mR^N)$,  
\begin{enumerate}
\item[$i)$] There exists $(u_n) \to u$ in $L^p(\mR^N)$ such that
$$
\limsup_{n \to + \infty} J_n(u_n) \le K_{N, p} \Phi(u)
$$

\item[$ii)$] For all $(u_n) \to u$ in $L^p(\mR^n)$, we have 
$$
\liminf_{n \to + \infty} J_n(u_n) \ge K_{N, p} \Phi(u). 
$$
\end{enumerate}
Assertion $i)$ follows from \Cref{thm-BBM} by taking $u_n = u$ for $n \ge 1$. 
\hmr{The proof of assertion $(ii)$, following the arguments in the work of Brezis and Nguyen \cite{BrNg-16}, where related results were also considered, can be carried out as follows.} 
Let $(\varphi_k)_{k \ge 1}$ be a smooth non-negative sequence of approximations to the identity such that $\supp \varphi_k \subset B_{1/k}$.  Since
$$
|a|^p - |b|^p \le C_p |a-b|^p, 
$$ 
it follows that, for $u, v \in L^p(\mR^N)$,  
\be
|J_n(u) - J_n(v)| \le C_{N, p} J_n(u-v). 
\ee
We derive that
\begin{multline}
|J_n(\phi_k*u) - J_n(\phi_k*u_n)| \le C_{N, p} J_n(\phi_k*u -\phi_k*u_n) \\[6pt]
\mathop{\le}^{\Cref{thm-BBM}} C_{N, p} \| \nabla \phi_k*u - \nabla \phi_k*u_n \|_{L^p(\mR^N)} \le  C_{N, p, k} \| u_n - u  \|_{L^p(\mR^N)}.  
\end{multline}
This implies 
\be
J_n(\phi_k*u) \le J_n(\phi_k* u_n) + C_{N, p, k} \| u_n - u  \|_{L^p(\mR^N)} \le J_n(u_n) + C_{N, p, k} \| u_n - u  \|_{L^p(\mR^N)}
\ee
By letting $n \to + \infty$, we obtain 
\be
K_{N, p} \Phi(\phi_k*u) \le \liminf_{n \to + \infty} J_n(u_n). 
\ee
By letting $k \to + \infty$, we obtain assertion $ii)$.

\medskip 
We next discuss the setting given in \Cref{thm-I} and \Cref{thm-B}. Set  
\be
\ka_{N, p, \gamma} : = \left\{ \begin{array}{cl} \inf \liminf_{\lambda \to 0_+}  \Phi_\lambda(h_\lambda)
& \mbox{ if } \gamma < 0,  \\[6pt]
\inf \liminf_{\lambda \to + \infty}  \Phi_\lambda(h_\lambda) & \mbox{ if } \gamma >0, 
\end{array} \right. 
\ee
where the infimum is taken over all families of measurable functions \( (h_\lambda)_{\lambda \in \mR} \) defined on the open unit cube $Q = (0, 1)^N$ such that \( h_\lambda \to h(x) \equiv \frac{x_1 + \dots + x_N}{\sqrt{N}} \) in Lebesgue measure on \( Q \) as \( \lambda \to 0 \) for $\gamma < 0$ and as $\lambda \to + \infty$ for $\gamma > 0$. It is clear that $\ka_{N, p, \gamma} \ge 0$. 

\medskip 
The following result was obtained in \cite{Ng-07,Ng-11}. 

\begin{theorem}[Nguyen] \label{thm-Gamma-I} Let $N \ge 1$, $p \ge 1$, and $\gamma = -p$. Then $0 < \ka_{N, p, \gamma} < \frac{K_{N, p}}{p}$ and   
\be
\Phi_\lambda \mbox{ $\Gamma$-converges to } \ka_{N, \gamma, p} \Phi \mbox{ in } L^p(\mR^N) \mbox{ as } \lambda \to 0.  
\ee
\end{theorem}

\begin{remark} \rm As a consequence of \Cref{thm-Gamma-I}, the $\Gamma$-limit is strictly smaller than the pointwise limit viewing \Cref{thm-I}. This fact is quite surprising. The proof of the fact $\ka_{N, p, -p} > 0$ is based on the work of Bourgain and Nguyen \cite{BN-06}, see also the proof of \Cref{thm-neg-1p} in \Cref{sect-thm-neg-1p}. The proof of the fact $\ka_{N, p, - p} < \frac{K_{N, p}}{p}$ is based on the construction of an example. The proof of the $\Gamma$-convergence is based on essentially the monotonicity property of $\Phi_\lambda$ in the case $\gamma = -p$. The analysis is delicate but quite robust and can be extended to a more general functionals, which are different from $\Phi_\delta$, by Brezis and Nguyen \cite{BN-18}. In the one dimensional case, we can obtain even a more general result \cite{BrNg-20}. 
\end{remark}

\begin{remark} \rm  The explicit value of $\ka_{N, p, -p}$ was conjectured in \cite{Ng-07} for $N=1$. This value was later confirmed by Antonucci, Gobbino, Migliorini, and Picenni  \cite{AGMP-20} where the case $N \ge 2$ was also established. 
\end{remark}

\begin{remark} \rm $\Gamma$-convergence results for functionals reminiscent of $\Phi_\lambda$ in the case $\gamma = - p$ were studied in by Ambrosio, De Philippis, and Martinazzi \cite{ADM-11}. 
\end{remark}

We conjecture that 
\begin{conjecture} \label{conj-Gamma} Let $N \ge 1$, $p \ge 1$, and $\gamma \in \mR \setminus \{0 \}$. Then 
\be
\Phi_\lambda \mbox{ $\Gamma$-converges to } \ka_{N, p, \gamma} \Phi \mbox{ in } L^p(\mR^N) \mbox{ for } \lambda \to + \infty \mbox{ if $\gamma >0$} \mbox{ and for } \lambda \to 0_+ \mbox{ if } \gamma < 0. 
\ee
\end{conjecture}

Concerning \Cref{conj-Gamma}, we have the following results. 

\begin{theorem} \label{thm-Gamma-Phi} Let $N \ge 1$, $p \ge 1$, and $-1 < \gamma < 0$. Then $\Phi_{\hmr{\lambda}}$ $\Gamma$-converges to $0$ in $L^p (\mR^N)$. 
\end{theorem}

\begin{proof} It suffices to prove that for each $u \in L^p(\mR^N)$, there exists $(u_\lambda)_{\lambda \in (0, 1)}$ such that $u_\lambda \to u$ in $L^p(\mR^N)$ and  $\Phi_\lambda(u_\lambda) \to 0$ as $\lambda \to 0_+$. For notational ease, we only consider the case $N=2$. The general case follows similarly. 

To this end, we first show this for each $u \in C^\infty_c(\mR^2)$ with compact support.  Let $m \in \mN$ be such that 
$$
\supp u \subset [-m, m]^2. 
$$
For $k \in \mN$, denote 
$$
\cQ_k = \Big\{\mbox{open cubes } Q_{i, j, k}: = ( i 2^{-k}, (i+1) 2^{-k}) \times (j2^{-k},  (j+1)2^{-k}) \subset \mR^2; i, j \in \mZ \Big\}.
$$

Set 
$$
u_{k} (x) = \sum_{-m 2^k \le i, j \le m 2^k} v_{i, j, k}(x),  
$$
where 
$$
v_{i, j, k} = u(x_{i, j,k}) \mathds{1}_{Q_{i, j, k}} \mbox{ with } x_{i, j, k} = \Big( (i+1/2) 2^{-k}, (j+1/2) 2^{-k} \Big). 
$$
We have, for some positive constants $c_{k}$ independent of $\hmr{\lambda}$,  
$$
\Phi_{c_{k}\lambda} (u_k) \le \sum_{-m 2^k \le i, j \le m 2^k} \Phi_{\lambda} (v_{i, j, k}). 
$$
Similar to the proof of \Cref{pro-CE}, we have 
$$
\lim_{\lambda \to 0_+} \Phi_{\lambda} (v_{i, j, k}) = 0. 
$$
We then derive that there exists a family $(u_\lambda)_{\lambda \in (0, 1)}$ such that $u_\lambda \to u$ in $L^p(\mR)$ and 
$$
\lim_{\lambda \to 0_+} \Phi_\lambda (u_\lambda) = 0. 
$$

Since for each $u \in L^p(\mR^N)$, there exist a sequence $(U_n) \subset C^\infty_c(\mR^N)$ such that $U_n \to u$ in $L^p(\mR^N)$.  The conclusion then follows from the case $u \in C^\infty_c(\mR^N)$ by a standard approximation. 
\end{proof}

\begin{remark}\rm
Brezis, Seeger, Van Schaftingen, and Yung \cite{BSVY-24} (see also \cite{Brezis-23}) asked whether the $\Gamma$-convergence holds in $L^1_{\loc}$ and the limit is $\Phi$ up to a positive constant. \Cref{thm-Gamma-Phi} gives a negative answer to this question in the case $-1 < \gamma < 0$. 
\end{remark}

\begin{remark} \rm  It would be interesting to \hmr{revisit} the arguments in \cite{Ng-11,BN-18} to establish \Cref{conj-Gamma}. 

\end{remark}

\section{Inequalities related to the nonlocal functionals} \label{sect-Inequality}

Since $\Phi_\lambda$ characterize Sobolev norms and the total variations, it is natural to ask whether or not one can obtain  properties of Sobolev spaces using the information of  $\Phi_\lambda$ instead of the one of $\Phi$. We addressed this question in the case $\gamma = -p$ \cite{Ng-Sob-11}. In particular, it was shown \cite{Ng-Sob-11} that a variant of Poincar\'e's inequality holds in this case. 

\begin{theorem}[Nguyen] \label{thm-Poincare}
Let $N \ge 1$, $p \ge 1$,  and $u$ be a real measurable function defined in a ball $B
\subset \mR^N$.  We have
\begin{equation}\label{Poincare}
\int_B \int_B |u(x) - u(y)|^p \, dx \, dy \le C_{N, p} \left(
|B|^{\frac{N+p}{N}}\mathop{\int_B \int_B}_{|\hmr{u}(x) - \hmr{u}(y)|
> \delta} \frac{\delta^p}{|x-y|^{N+p}} \, dx \, dy + \delta^p
|B|^{2} \right). 
\end{equation}
\end{theorem}

The proof of \Cref{thm-Poincare} is based on the arguments of Bourgain and Nguyen \cite{BN-06} used in the proof of \Cref{lem-fund}, and the inequalities associated with $BMO$-functions due to John and Nirenberg \cite{JN-61}.  

\medskip 
Applying Theorem~\ref{thm-Poincare}, one can derive that  $u \in BMO(\mR^N)$, the space of all functions of bounded mean oscillation defined in $\mR^N$ if $u \in L^1(\mR^N)$  and $\Phi_\delta(u) < +\infty$  for $\gamma = -p$ and $p =N$,  and for some $\lambda >0$.
Moreover, there exists a positive constant $C$, depending only on $N$, such that, for $\gamma = - p$ and $p = N$, 
\begin{equation*}
|g|_{BMO} : = \sup_{B} \fint_B \fint_B |g(x) - g(y)| \, dx \, dy
\le C \left( \Phi_{\lambda}^\frac{1}{N} (u) + \lambda \right),
\end{equation*}
where the supremum is taken over all balls of $\mR^N$. In a joint work with Brezis \cite{BN-VMO-11}, we also show that if $u \in L^1(\mR^N)$ and $\Phi_\lambda (u) < + \infty$ with $\gamma =- p$ and  $p=N$ for all $\lambda >0$, then $g \in VMO(\mR^N)$, the spaces of all functions of vanishing mean oscillation. 
\medskip

Using \Cref{thm-Poincare}, we can establish variants of Sobolev's inequalities and Rellich-Kondrachov's compactness criterion \cite{Ng-Sob-11}. The proof of Rellich-Kondrachov's compactness criterion using \Cref{thm-Poincare} is quite standard. The idea of the proof of the Sobolev inequalities using the Poincar\'e inequalities is as follows. We first establish the Sobolev inequalities for the weak type from the Poincar\'e inequalities using covering lemmas. We then used the truncation arguments due to \hmr{Maz'ya}, see, e.g., \cite{MazyaSobolev}, and an inequality related to sharp maximal functions due to Fefferman and Stein, see, e.g., \cite{Stein-70}, to obtain the desired estimates from the weak-type ones. This kind \hmr{of} arguments have been extended to obtain the full range of Gagliardo \& Nirenberg and Caffarelli \& Kohn \& Nirenberg interpolation inequalities associated with Coulomb-Sobolev spaces \cite{Ng-Mallick2}, a result obtained in a collaboration with Mallick.

In another direction, one can also derive variants of the Hardy inequalities and the Caffarelli \& Kohn \& Nirenberg inequalities using the information of $\Phi_\lambda$ instead of $\Phi$ in the case $\gamma = -p$. This is given in a joint work with Squassina \cite{NS-Hardy-19}.  Interestingly, our proofs are quite elementary and mainly based on the Poincar\'e and Sobolev inequalities for an annulus; the integration-by-part arguments are not required. Our analysis is inspired from the harmonic one, nevertheless, instead of using dyadic decomposition for the frequency, we do it for the space variables. These arguments have also used by us to obtain the full range of the  Caffarelli \& Kohn \& Nirenberg inequalities for fractional Sobolev's spaces \cite{NS-18}, which generalizes the Caffarelli \& Kohn \& Nirenberg inequalities in \cite{CKN}.

\medskip
In this direction, concerning $\Phi_\lambda$,  we ask the following question viewing \Cref{pro-CE}.  

\begin{question} Let $Q$ be a unit cube of $\mR^N$, and let $p \ge 1$, $\gamma \in \mR \setminus (-1, 0]$ be such that $\ka_{N, p, \gamma} > 0$. Is it true that 
$$
\iint_{Q\times Q} |u(x) - u(y)|^p \, dx \, dy \le C_{N, p} \Big( \Phi_\lambda (u) + \lambda^p \Big) \mbox{ for } u \in L^p(Q)? 
$$
\end{question}

{\color{black}
\noindent {\bf Acknowledgement:} The author is grateful to the referee for a careful reading of the manuscript and for valuable suggestions. }

\providecommand{\bysame}{\leavevmode\hbox to3em{\hrulefill}\thinspace}
\providecommand{\MR}{\relax\ifhmode\unskip\space\fi MR }
\providecommand{\MRhref}[2]{%
  \href{http://www.ams.org/mathscinet-getitem?mr=#1}{#2}
}
\providecommand{\href}[2]{#2}


\begin{thebibliography}{10}

\bibitem{Adams-75}
Robert~A. Adams, \emph{Sobolev spaces}, Pure and Applied Mathematics, vol. Vol.
  65, Academic Press [Harcourt Brace Jovanovich, Publishers], New York-London,
  1975. \MR{450957}

\bibitem{ADM-11}
Luigi Ambrosio, Guido De~Philippis, and Luca Martinazzi,
  \emph{Gamma-convergence of nonlocal perimeter functionals}, Manuscripta Math.
  \textbf{134} (2011), no.~3-4, 377--403. \MR{2765717}

\bibitem{AGMP-20}
Clara Antonucci, Massimo Gobbino, Matteo Migliorini, and Nicola Picenni,
  \emph{Optimal constants for a nonlocal approximation of {S}obolev norms and
  total variation}, Anal. PDE \textbf{13} (2020), no.~2, 595--625. \MR{4078237}

\bibitem{BBM-94}
Fabrice Bethuel, Ha\"im Brezis, and Fr\'ed\'eric H\'elein,
  \emph{Ginzburg-{L}andau vortices}, Progress in Nonlinear Differential
  Equations and their Applications, vol.~13, Birkh\"auser Boston, Inc., Boston,
  MA, 1994. \MR{1269538}

\bibitem{BBM-01}
Jean Bourgain, Ha\"im Brezis, and Petru Mironescu, \emph{Another look at
  {S}obolev spaces}, Optimal control and partial differential equations, IOS,
  Amsterdam, 2001, pp.~439--455. \MR{3586796}

\bibitem{BBM-04}
Jean Bourgain, Ha\"im Brezis, and Petru Mironescu, \emph{{Complements to the
  paper:``Lifting, Degree, and Distributional Jacobian Revisited"}},  (2004),
  unpublished, {http://hal.archives-ouvertes.fr/docs/00/74/76/68/PDF}.

\bibitem{BBM-05}
Jean Bourgain, Ha\"im Brezis, and Petru Mironescu, \emph{Lifting, degree, and
  distributional {J}acobian revisited}, Comm. Pure Appl. Math. \textbf{58}
  (2005), no.~4, 529--551. \MR{2119868}

\bibitem{BBN-05}
Jean Bourgain, Ha\"im Brezis, and Hoai-Minh Nguyen, \emph{A new estimate for
  the topological degree}, C. R. Math. Acad. Sci. Paris \textbf{340} (2005),
  no.~11, 787--791. \MR{2139888}
  
\bibitem{BN-06}
Jean Bourgain and Hoai-Minh Nguyen, \emph{A new characterization of {S}obolev
  spaces}, C. R. Math. Acad. Sci. Paris \textbf{343} (2006), no.~2, 75--80.
  \MR{2242035}

\bibitem{Brezis-FA}
Ha\"{\i}m Brezis, \emph{Functional analysis, {S}obolev spaces and partial
  differential equations}, Universitext, Springer, New York, 2011. \MR{2759829}
  
  \bibitem{Brezis-02}
Ha\"im Brezis, \emph{How to recognize constant functions. {A} connection with
  {S}obolev spaces}, Uspekhi Mat. Nauk \textbf{57} (2002), no.~4(346), 59--74.
  \MR{1942116}
  
  \bibitem{Brezis-23}
Ha\"im Brezis, \emph{Some of my favorite open problems}, Atti Accad. Naz.
  Lincei Rend. Lincei Mat. Appl. \textbf{34} (2023), no.~2, 307--335.
  \MR{4668549}
  
\bibitem{BN-VMO-11}
Ha\"{\i}m Brezis and Hoai-Minh Nguyen, \emph{On a new class of functions
  related to {VMO}}, C. R. Math. Acad. Sci. Paris \textbf{349} (2011), no.~3-4,
  157--160. \MR{2769899}

\bibitem{BrNg-16}
Ha\"im Brezis and Hoai-Minh Nguyen, \emph{The {BBM} formula revisited}, Atti
  Accad. Naz. Lincei Rend. Lincei Mat. Appl. \textbf{27} (2016), no.~4,
  515--533. \MR{3556344}

\bibitem{BN-18}
Ha\"im Brezis and Hoai-Minh Nguyen, \emph{Non-local functionals related to
  the total variation and connections with image processing}, Ann. PDE
  \textbf{4} (2018), no.~1, Paper No. 9, 77. \MR{3749763}

\bibitem{BrNg-20}
Ha\"im Brezis and Hoai-Minh Nguyen, \emph{{$\Gamma$}-convergence of non-local, non-convex functionals in
  one dimension}, Commun. Contemp. Math. \textbf{22} (2020), no.~7, 1950077,
  27. \MR{4135011}
  
\bibitem{BSVY-222}
Ha\"im Brezis, Andreas Seeger, Jean Van~Schaftingen, and Po-Lam Yung,
  \emph{Sobolev spaces revisited}, Atti Accad. Naz. Lincei Rend. Lincei Mat.
  Appl. \textbf{33} (2022), no.~2, 413--437. \MR{4482042}

\bibitem{BSVY-24}
Ha\"im Brezis, Andreas Seeger, Jean Van~Schaftingen, and Po-Lam Yung,
  \emph{Families of functionals representing {S}obolev norms}, Anal. PDE
  \textbf{17} (2024), no.~3, 943--979. \MR{4736522}

\bibitem{BVY-21b}
Ha\"im Brezis, Jean Van~Schaftingen, and Po-Lam Yung, \emph{Going to {L}orentz
  when fractional {S}obolev, {G}agliardo and {N}irenberg estimates fail}, Calc.
  Var. Partial Differential Equations \textbf{60} (2021), no.~4, Paper No. 129,
  12. \MR{4279388}

\bibitem{BVY-21}
Ha\"im Brezis, Jean Van~Schaftingen, and Po-Lam Yung, \emph{A surprising formula for {S}obolev norms}, Proc. Natl. Acad.
  Sci. USA \textbf{118} (2021), no.~8, Paper No. e2025254118, 6. \MR{4275122}



\bibitem{CKN}
Luis Caffarelli, Robert Kohn, and Louis Nirenberg, \emph{First order interpolation
  inequalities with weights}, Compositio Math. \textbf{53} (1984), no.~3,
  259--275. \MR{768824}

\bibitem{Davila-02}
Juan D\'{a}vila, \emph{On an open question about functions of bounded
  variation}, Calc. Var. Partial Differential Equations \textbf{15} (2002),
  no.~4, 519--527. \MR{1942130}

\bibitem{DM-22}
{\'O}scar Dom{\'\i}nguez and Mario Milman, \emph{New Brezis--Van
  Schaftingen--Yung--Sobolev type inequalities connected with maximal
  inequalities and one parameter families of operators}, Advances in
  Mathematics \textbf{411} (2022), 108774.


\bibitem{EGMeasure}
Lawrence~C. Evans and Ronald~F. Gariepy, \emph{Measure theory and fine
  properties of functions}, Studies in Advanced Mathematics, CRC Press, Boca
  Raton, FL, 1992. \MR{1158660}
  
  
\bibitem{GB-25}
Massimo Gobbino and Nicola Picenni, \emph{Gamma-liminf estimate for a class of
  non-local approximations of {S}obolev and {BV} norms}, J. Funct. Anal.
  \textbf{289} (2025), no.~9, Paper No. 111106. \MR{4925908}

\bibitem{JN-61}
Fritz John and Louis Nirenberg, \emph{On functions of bounded mean
  oscillation}, Comm. Pure Appl. Math. \textbf{14} (1961), 415--426.
  \MR{131498}
  
\bibitem{LS-11}
Giovanni Leoni and Daniel Spector, \emph{Characterization of {S}obolev and
  {$BV$} spaces}, J. Funct. Anal. \textbf{261} (2011), no.~10, 2926--2958.
  \MR{2832587}

\bibitem{Ng-Mallick2}
Arka Mallick and Hoai-Minh Nguyen, \emph{Gagliardo-{N}irenberg and
  {C}affarelli-{K}ohn-{N}irenberg interpolation inequalities associated with
  {C}oulomb-{S}obolev spaces}, J. Funct. Anal. \textbf{283} (2022), no.~10,
  Paper No. 109662, 33. \MR{4474840}

\bibitem{MazyaSobolev}
Vladimir Maz'ya, \emph{Sobolev spaces with applications to elliptic partial
  differential equations}, augmented ed., Grundlehren der mathematischen
  Wissenschaften [Fundamental Principles of Mathematical Sciences], vol. 342,
  Springer, Heidelberg, 2011. \MR{2777530}

\bibitem{MS-02}
Vladimir Mazya and Tatyana~O. Shaposhnikova, \emph{{On the Bourgain,
  Brezis, and Mironescu theorem concerning limiting embeddings of fractional
  Sobolev spaces}}, J. Funct. Anal. \textbf{195} (2002), no.~2, 230--238.
  \MR{1940355}

\bibitem{NgSob1}
Hoai-Minh Nguyen, \emph{Some new characterizations of {S}obolev spaces}, J.
  Funct. Anal. \textbf{237} (2006), no.~2, 689--720. \MR{2230356}

\bibitem{Ng-07}
Hoai-Minh Nguyen, \emph{{$\Gamma$}-convergence and {S}obolev norms}, C. R. Math. Acad.
  Sci. Paris \textbf{345} (2007), no.~12, 679--684. \MR{2376638}


\bibitem{Ng-Opt}
Hoai-Minh Nguyen, \emph{Optimal constant in a new estimate for the degree}, J. Anal.
  Math. \textbf{101} (2007), 367--395. \MR{2346551}
  
\bibitem{NgSob2}
Hoai-Minh Nguyen, \emph{Further characterizations of {S}obolev spaces}, J. Eur. Math.
  Soc. (JEMS) \textbf{10} (2008), no.~1, 191--229. \MR{2349901}

\bibitem{Ng-11}
Hoai-Minh Nguyen, \emph{{$\Gamma$}-convergence, {S}obolev norms, and {BV} functions},
  Duke Math. J. \textbf{157} (2011), no.~3, 495--533. \MR{2785828}

\bibitem{Ng-Sob-11}
Hoai-Minh Nguyen, \emph{Some inequalities related to {S}obolev norms}, Calc. Var.
  Partial Differential Equations \textbf{41} (2011), no.~3-4, 483--509.
  \MR{2796241}
  
  \bibitem{Ng-refine}
Hoai-Minh Nguyen, \emph{A refined estimate for the topological degree}, C. R. Math.
  Acad. Sci. Paris \textbf{355} (2017), no.~10, 1046--1049. \MR{3716480}

\bibitem{NPSV-18}
Hoai-Minh Nguyen, Andrea Pinamonti, Marco Squassina, and Eugenio Vecchi,
  \emph{New characterizations of magnetic {S}obolev spaces}, Adv. Nonlinear
  Anal. \textbf{7} (2018), no.~2, 227--245. \MR{3794886}

\bibitem{NS-18}
Hoai-Minh Nguyen and Marco Squassina, \emph{Fractional
  {C}affarelli-{K}ohn-{N}irenberg inequalities}, J. Funct. Anal. \textbf{274}
  (2018), no.~9, 2661--2672. \MR{3771839}

\bibitem{NS-19}
Hoai-Minh Nguyen and Marco Squassina, \emph{On anisotropic {S}obolev spaces}, Commun. Contemp. Math.
  \textbf{21} (2019), no.~1, 1850017, 13. \MR{3904643}

\bibitem{NS-Hardy-19}
Hoai-Minh Nguyen and Marco Squassina, \emph{On {H}ardy and {C}affarelli-{K}ohn-{N}irenberg inequalities}, J.
  Anal. Math. \textbf{139} (2019), no.~2, 773--797. \MR{4041120}

\bibitem{Picenni-24}
Nicola Picenni, \emph{New estimates for a class of non-local approximations of
  the total variation}, J. Funct. Anal. \textbf{287} (2024), no.~1, Paper No.
  110419, 21. \MR{4733863}

\bibitem{PSV-19}
Andrea Pinamonti, Marco Squassina, and Eugenio Vecchi, \emph{Magnetic
  {BV}-functions and the {B}ourgain-{B}rezis-{M}ironescu formula}, Adv. Calc.
  Var. \textbf{12} (2019), no.~3, 225--252. \MR{3975602}

\bibitem{Poliakovsky-22}
Arkady Poliakovsky, \emph{Some remarks on a formula for {S}obolev norms due to
  {B}rezis, {V}an {S}chaftingen and {Y}ung}, J. Funct. Anal. \textbf{282}
  (2022), no.~3, Paper No. 109312, 47. \MR{4339013}

\bibitem{Ponce-04}
Augusto~C. Ponce, \emph{A new approach to {S}obolev spaces and connections to
  {$\Gamma$}-convergence}, Calc. Var. Partial Differential Equations
  \textbf{19} (2004), no.~3, 229--255. \MR{2033060}

\bibitem{PS-17}
Augusto~C. Ponce and Daniel Spector, \emph{On formulae decoupling the total
  variation of {BV} functions}, Nonlinear Anal. \textbf{154} (2017), 241--257.
  \MR{3614653}

\bibitem{Stein-70}
Elias~M. Stein, \emph{Singular integrals and differentiability properties of
  functions}, Princeton Mathematical Series, vol. No. 30, Princeton University
  Press, Princeton, NJ, 1970. \MR{290095}

\end{thebibliography}
\end{document}